\documentclass{gtart}
\usepackage{amsmath,amssymb}
\let\Bbb\mathbb

\input gtoutput
\volumenumber{3}\papernumber{8}\volumeyear{1999}
\pagenumbers{167}{210}\published{4 July 1999}
\proposed{Robion Kirby}
\seconded{Gang Tian, Tomasz Mrowka}
\received{26 July 1998}
\accepted{8 May 1999}

\renewcommand{\theequation}{\thesection.\arabic{equation}}

\newcommand{\eqnum}{\addtocounter{equation}{1}%
\hspace\fill{\rm (\theequation)}}

\newtheorem{theorem}{Theorem}[section]
\newtheorem{lemma}{Lemma}[section]

\newtheorem{proposition}{Proposition}[section]

\newcommand{\x}{\times}
\newcommand{\<}{\langle}
\renewcommand{\>}{\rangle}
\newcommand{\R}{{\Bbb R}}
\newcommand{\Z}{{\Bbb Z}}

\newcommand{\ot}{\otimes}

\newcommand{\C}{{\Bbb C}}
\newcommand{\dt}{\cdot}
\newcommand{\bs}{\medskip}

\renewcommand{\a}{\alpha}
\renewcommand{\b}{\beta}
\renewcommand{\d}{\delta}
\newcommand{\D}{\Delta}
\newcommand{\e}{\varepsilon}
\newcommand{\g}{\gamma}

\renewcommand{\l}{\lambda}
\renewcommand{\L}{\Lambda}

\renewcommand{\phi}{\varphi}
\newcommand{\s}{\sigma}

\renewcommand{\O}{\Omega}
\renewcommand{\o}{\omega}
\newcommand{\z}{\zeta}
\renewcommand{\i}{\infty}
\newcommand{\p}{\partial}

\title{Seiberg--Witten Invariants and Pseudo-Holomorphic\\Subvarieties
for Self-Dual, Harmonic 2--Forms}
\shorttitle{Seiberg--Witten Invariants and Pseudo-Holomorphic Subvarieties}

\author{Clifford Henry Taubes}

\address{Department of Mathematics \\
Harvard University \\
Cambridge, MA 02138, USA}

\email{chtaubes@abel.math.harvard.edu}

\begin{abstract}
 A smooth, compact 4--manifold with a Riemannian metric and $b^{2+}\ge 1$
has a non-trivial, closed, self-dual 2--form.  If the metric is
generic, then the zero set of this form is a disjoint union of
circles.  On the complement of this zero set, the symplectic form and
the metric define an almost complex structure; and the latter can be
used to define pseudo-holomorphic submanifolds and subvarieties.  The
main theorem in this paper asserts that if the 4--manifold has a non
zero Seiberg--Witten invariant, then the zero set of any given
self-dual harmonic 2--form is the boundary of a pseudo-holomorphic 
subvariety in its complement.
\end{abstract}

\asciiabstract{A smooth, compact 4-manifold with a Riemannian metric 
and b^(2+) > 0 has a non-trivial, closed, self-dual 2-form.  If the
metric is generic, then the zero set of this form is a disjoint union
of circles.  On the complement of this zero set, the symplectic form
and the metric define an almost complex structure; and the latter can
be used to define pseudo-holomorphic submanifolds and subvarieties.
The main theorem in this paper asserts that if the 4-manifold has a
non zero Seiberg-Witten invariant, then the zero set of any given
self-dual harmonic 2-form is the boundary of a pseudo-holomorphic
subvariety in its complement.}




\begin{document}

\maketitlepage

\section{Introduction}
\label{sec1}

Let $X$ be a compact, oriented 4--dimensional manifold with Betti number $b^{2+}
\geq 1$. Choose a Riemannian metric, $g$, for $X$ and Hodge theory provides
a $b^{2+}$--dimensional space of self-dual, harmonic 2--forms. Let $\o$ be
such a self-dual harmonic 2--form. At points where $\o\neq 0$, the endomorphism
$J=\sqrt{2} |\o|^{-1}g^{-1}\o$ of $TX$ has square equal to minus the
identity and thus defines an almost complex structure. The latter can be
used to define, after Gromov \cite{Gr}, the notion of a pseudo-holomorphic
curve in the
complement of the zero set of $\o$. This last notion can be generalized with
the following definition.

\bs{\bf Definition 1.1}\qua Let $Z\subset X$ denote the zero set of $\o$. A
subset $C\subset X-Z$ will be called {\it finite energy, pseudo-holomorphic
subvariety}
when the following requirements are met:
\begin{itemize}
\item There is a complex curve $C_0$ (not necessarily compact or connected)
together with a proper, pseudo-holomorphic map $\phi\co  C_0\to X-Z$ such that
$\phi(C_0)=C$.
\item There is a countable set $\L_0\subset C_0$ which has no accumulation
points and is such that $\phi$ embeds $C_0-\L_0$.
\item The integral of $\phi^*\o$ over $C_0$ is finite.
\end{itemize}

When $Z=\emptyset$ and so the form $\o$ is symplectic, then the main
theorem in \cite{T1} asserts that pseudo-holomorphic subvarieties exist
when the
Seiberg--Witten
invariants of $X$ are not all zero. (A different sort of existence theorem
for pseudo-holomorphic curves has been given by Donaldson \cite{D}.)

The purpose of this paper is to provide a generalization of the existence
theorem in \cite{T1} to the case where $Z \neq \emptyset$. The statement of this
generalization is cleanest in the case where $\o$ vanishes transversely.
This turns out to be the generic situation, see eg~\cite{LeB},\cite{Ho}.
In this
case $Z$ is a union of embedded circles. The following theorem summarizes
the existence theorem in this case.

\setcounter{theorem}{1}
\begin{theorem}
\label{th1.2}
Let $X$ be a compact, oriented, Riemannian
4--manifold with $b^{2+}\geq 1$ and with a non-zero Seiberg--Witten
invariant. Let $\o$ be a
self-dual, harmonic 2--form which vanishes transversely. Then there is a
finite energy,
pseudo-holomorphic subvariety $C\subset X-Z$ with the property that $C$ has
intersection
number 1 with every linking 2--sphere of $Z$.
\end{theorem}

Some comments are in order. In the case where $b^{2+}=1$, the
Seiberg--Witten invariants require
a choice of ``chamber" for their definition. Implicit in the statement of
Theorem $1.1$ is that the chamber in question is defined by a perturbation of
the equations which is constructed from the chosen form $\o$. This chamber
is described in more detail in the next section.

Here is a second comment: With regard to the sign of the intersection
number between $C$ and the linking 2--spheres of $Z$, remark that a
pseudo-holomorphic
subvariety is canonically oriented at its smooth points by the restriction
of $\o$. Meanwhile, the linking 2--spheres of $Z$ are oriented as follows:
First, use the assumed orientation of $TX$ to orient the bundle of self-dual
2--forms. Second, use the differential of $\o$ along $Z$ to identify the
normal bundle of $Z$ with this same bundle of 2--forms. This orients the
normal bundle to $Z$ and thus $TZ$.

Here is the final comment: Away from $Z$, the usual regularity theorems for
pseudo-holomorphic curves (as in \cite{MS},\cite{Pa},\cite{PW} or
\cite{Ye}) describe the
structure of a
finite energy, pseudo-holomorphic subvariety. Basically, such a subvariety
is no more
singular than an algebraic curve in $\C^2$. However, since the almost
complex structure $J$ is singular along $Z$, there are serious questions
about the regularity near $Z$ of a finite energy, pseudo-holomorphic
subvariety. In this
regard, \cite{T2} provides a first step towards describing the general
structure.
In \cite{T2} the metric is restricted near $Z$ to have an especially simple
form
and for this restricted metric the story, as developed to date, is as
follows: All but finitely many points on $Z$ have a ball neighborhood which
the finite energy subvariety intersects in a finite number of disjoint
components. Moreover, the closure of each such component in this ball is a
smoothly embedded, closed half-disc whose straight edge coincides with $Z$.
(There are no obstructions to realizing the special metrics of \cite{T2}.)

The next section also provides some examples where the subvarieties of
Theorem \ref{th1.2} are easy to see.

Note that Theorem 2.3, below, gives an existence assertion without assuming
the transversality of $\o$ at zero. Also, Theorem 2.2 below is a stronger
version of Theorem \ref{th1.2}.

This introduction ends with an open problem for the reader (see also
\cite{T3}).

\bs{\bf Problem}\qua The proof of Theorem \ref{th1.2} suggests that the
Seiberg--Witten invariants
of any
$b^{2+}$ positive 4--manifold can be computed via a creative algebraic count
of the finite energy, pseudo-holomorphic subvarieties which (homologically)
bound the zero
set of a non-trivial, self-dual harmonic form. Such a count is known in the
case where the form is nowhere vanishing (see \cite{T4},\cite{T5}). 
Another case which
is well understood has the product metric on $X=S^1\x M$, where $M$ is a
compact, oriented, Riemannian 3--manifold with positive first Betti number
(\cite{HL1},\cite{HL2},\cite{Tu}). The problem is to find such a count which applies
for any
compact, $b^{2+}$ positive 4--manifold.

\section{Basics}
\label{sec2}

The purpose of this section is to review some of the necessary
background for the Seiberg--Witten equations and for the study of
pseudo-holomorphic
subvarieties.

\bs{\bf a) The Seiberg--Witten equations}

The Seiberg--Witten equations are discussed in numerous sources to date,
so  the discussion here will be brief. The novice can consult the book by
Morgan \cite{Mor} or the forthcoming book by Kotschick, Kronheimer and Mrowka
\cite{KKM}.

To begin, suppose for the time being that $X$ is an oriented, Riemannian
4--manifold. The chosen metric on $X$ defines the principal SO(4) bundle
$Fr\to X$ of oriented, orthonormal frames in $TX$. A {\it spin${}^{\C}$
structure} on $X$ is a lift (or, more properly, an equivalence class of
lifts) of $Fr$ to a principal Spin${}^{\C}$(4) bundle $F\to X$. In this
regard the reader should note the identifications
\setcounter{equation}{0}
\begin{itemize}
\item ${\mbox{Spin}}^{\C}(4) = (SU(2)\x SU(2)\x U(1))/\{\pm 1\}$,
\item $SO(4)  = (SU(2)\x SU(2))/\{\pm 1\}$ \eqnum
\end{itemize}
%\begin{equation}
%\label{eq2.1}
%\end{equation}
with the evident group homomorphism from the former to the latter which
forgets the factor of $U(1)$ in the top line above. Remark that there exist,
in any event Spin${}^{\C}$ structures on 4--manifolds.  Moreover, the set
$\cal S$ of Spin${}^{\C}$ structures is naturally metric independent and has
the structure of a principal homogeneous space for the additive group
$H^2(X;\Z)$.

With the preceding understood, fix a Spin${}^{\C}$ structure $F\to X$. Then
$F$ can be used to construct three useful associated vector bundles, $S_+$,
$S_-$ and $L$. The first two are associated via the representations
$s_{\pm}\co {\mbox{Spin}}^{\C}(4)\to U(2)= (SU(2)\x U(1))/\{\pm 1\}$ which forgets
one or the other factor of $SU(2)$ in the top line of (2.1). Thus
$S_{\pm}$ are $\C^2$ vector bundles over $X$ with Hermitian metrics.
Meanwhile,
$L={\mbox{det}}(S_+)={\mbox{det}}(S_-)$ is associated to $F$ via the
representation of ${\mbox{Spin}}^{\C}(4)$ on $U(1)$ which forgets both
factors of $SU(2)$ in the first line of (2.1). (By way of
comparison, the
$\R^3$ bundles $\L_+,\L_-\to X$ of self-dual and anti-self-dual 2--forms are
associated to $Fr$ via the representations of $SO(4)$ to $SO(3)$ which forget
one of the other factors of $SU(2)$ in the second line of (2.1).)

Note that the bundle $S_+\oplus S_-$ is a module for the Clifford algebra of
$TX$ in the sense that there is an epimorphism
$cl\co TX\to{\mbox{Hom}}(S_+,S_-)$ which obeys $cl^{\dagger}cl=-1$. The latter
will be thought of equally as a homomorphism from $S_+\ot TX$ to $S_-$. Note
that this homomorphism induces one, $cl_+$, from $\L_+$ to End$(S_+)$.

Now consider that the Seiberg--Witten equations constitute a system of
differential equations for a pair $(A,\psi)$, where $A$ is a hermitian
connection on the complex line bundle $L$ and where $\psi$ is a section of
$S_+$. These equations read, schematically:
\begin{equation}
\label{eq2.2}
\begin{array}{rll}
D_A\psi &=& 0 \\
F_A^+ & =& q(\psi)+\mu .
\end{array}
\end{equation}
In the first line above, $D_A$ is the Dirac operator as defined using the
connection $A$ and the Levi-Civita connection on $TX$.
Indeed, these two connections define a unique connection on $S_+$ and thus
a covariant derivative, $\nabla_A$, which takes a section of $S_+$ and
returns one of $S_+\ot T^*X$. With this understood, then $D_A\psi$ sends the
section $\psi$ of $S_+$ to the section cl$(\nabla_A\psi)$ of $S_-$.
In the second line of (\ref{eq2.2}), $F_A$ is the curvature 2--form of the
connection
$A$ on
$L$, this being an imaginary valued 2--form. Then $F^+_A$ is the projection
of $F_A$ onto $\L_+$. Meanwhile, $q(\dt)$ is the quadratic map from $S_+$ to
$i\dt\L_+$ which, up to a constant factor, sends $\eta\in S_+$ to the image
of $\eta\otimes\eta^{\dagger}$ under the adjoint of cl${}_+$. To be more
explicit about $q$, let $\{e^\nu\}_{1\leq\nu\leq 4}$ be an oriented,
orthonormal frame for $T^*X$. Then
$q(\eta)=-8^{-1}\sum_{\nu,\l}\<\eta,{\mbox{cl}}(e^{\nu})
{\mbox{cl}}(e^{\l})\eta\>e^{\nu}\wedge e^{\l}$ where $\< \ , \ \>$ denotes
the Hermitian inner product on $S_+$. Finally, $\mu$ in the second line of
(\ref{eq2.2}) is a favorite, imaginary valued, self-dual 2--form. (A different
choice for
$\mu$, as with a different choice for the Riemannian metric, will give a
different set of equations.)

\bs{\bf b) The Seiberg--Witten invariants}

Let $Q$ denote the cup product pairing of $H^2(X;\R)$ and let
$H^{2+}\subset H^2(X;\R)$ denote a maximum subspace on which $Q$ is positive
definite. Set $b^{2+}={\mbox{dim}}(H^{2+})$. Fix an orientation for the real
line $\L^{\mbox{\scriptsize}{top}}H^1(X;\R)
\otimes\L^{\mbox{\scriptsize}{top}}H^{2+}$. If $b^{2+} >1$,
then the Seiberg--Witten invariants as presented in \cite{W} constitute a
diffeomorphism
invariant map
$SW\co \cal S\to\Z$. Moreover, there is a straightforward generalization (see
\cite{T5}) which extends this invariant to
\begin{equation}
\label{eq2.3}
SW\co {\cal S} \to\L^*H^1(X;\Z) .
\end{equation}
In the case where $b^{2+} =1$, there is a diffeomorphism invariant as in
(\ref{eq2.3}) after the additional choice of an orientation for the line
$H^{2+}$.

In all cases, the map in (\ref{eq2.3}) is defined via a creative, algebraic
count  of the solutions of (\ref{eq2.2}). However, the particulars of the
definition of
$SW$ are not relevant to the discussion in this article except for the
following two facts:
\begin{itemize}
\item If $SW(s)\neq 0$, then there exists, for each choice of
metric $g$ and perturbing form $\mu$, at least one solution to (\ref{eq2.2})
as defined by the Spin${}^{\C}$ structure $s$.
\item In the case
where $b^{2+} =1$, the orientation of $H^{2+}$ defines a unique self-dual
harmonic 2--form $\o$ up to multiplication by the positive real numbers.
With this
understood, note that $SW$ in (\ref{eq2.3}) is computed by counting solutions
to (\ref{eq2.2}) in the special case where the perturbation $\mu$ in
(\ref{eq2.2}) has the form
$\mu=-i\dt r/4\dt\o+\mu_0$, where $\mu_0$ is a fixed, imaginary valued
2--form and where $r$ is taken to be very large. That is, the algebraic count
of solutions to (\ref{eq2.2}) stabilizes as $r$ tends to $+\i$, and the large
$r$ count is defined to be $SW$. \eqnum
\end{itemize}
%\begin{equation}
%\label{eq2.4}
%\end{equation}

\bs{\bf c) Near the zero set of a self-dual harmonic form}

Let $X$ be a compact, oriented, Riemannian 4--manifold with $b^{2+}>0$ and
suppose that $\o$ is a self-dual harmonic 2--form which vanishes
transversely. The
purpose of this subsection is to describe the local geometry of the zero set
$Z\equiv\o^{-1}(0)$, of the form $\o$.

To begin, note first that the non-degeneracy condition implies that $Z$ is
a union of embedded circles. Moreover, the transversal vanishing of $\o$
implies that its covariant derivative, $\nabla\o$, identifies the normal
bundle $N\to Z$ of $Z$ with the bundle $\L_+|_Z$ of self-dual 2--forms. As
$\L_+$ is oriented by the orientation of $X$, the homomorphism $\nabla\o$
orients $N$ with the declaration that it be orientation reversing. This
orientation of $N$ induces one on $Z$ if one adopts the convention that
$TX=TZ\oplus N$ (as opposed to $N\oplus TZ$).

With $Z$ now oriented, define $\tau\co \L_+\to N^*$ by the rule $\tau(u)\equiv
u(\p_0,\dt)$, where the $\p_0$ is the unit length oriented tangent vector to
$Z$. Note that $\tau$ is also an isomorphism. Moreover, the composition
$\tau \cdot \nabla\o\co N\to N^*$ defines a bilinear form on $N$ with negative
determinant. And, as $d\o=0$, this form is symmetric with trace zero and
thus $\tau\dt\nabla\o$ has everywhere three real eigenvalues, where two are
positive and one is negative. (Note here that $N$ inherits a fiber metric
with its identification as the orthogonal complement to $TZ$ in $TX|_Z$.)

Let $N_1\subset N$ denote the one-dimensional eigenbundle for
$\tau\dt\nabla\o$ which corresponds to the negative eigenvalue. Then use
$N_2\subset N$ to denote its orthogonal complement. With regard to $N_1$,
note that this bundle can be either oriented or not. Gompf has pointed out
that $1+b^2-b^1$ and the number of components of $Z$ for which
$N_1$ is oriented are equal modulo 2.

With the Riemannian geometry near $Z$ understood, consider now the almost
complex geometry in a neighborhood of $Z$. Here the almost complex structure
on $X-Z$ is defined by the endomorphism $J\equiv\sqrt{2}g^{-1}\o/|\o|$ with
$\o$ viewed as a skew symmetric homomorphism from $TX$ to $T^*X$ and with
$g^{-1}$ viewed as a symmetric homomorphism which goes the other way. As $J$
has square equal to minus the identity, $J$ decomposes
$T^*X|_{X-Z} \otimes \C  =T^{1,0}\oplus T^{0,1}$, where $T^{1,0}$ are the holomorphic
1--forms and $T^{0,1}$ the anti-holomorphic forms. The {\it canonical line
bundle}\,, $K$, for the almost complex structure is $\L^2T^{1,0}$.

Note that $J$ does not extend over $Z$. This failure is implied by the
following lemma:

\begin{lemma}
\label{lem2.1}
Let $\s\co N\to X$ denote the metric's exponential map
and let $N^0\subset N$ be an open ball neighborhood of the zero section
which is embedded by $\s$. Use $\s$ to identify $N^0$ with a neighborhood of
$Z$ in $X$. Let $p\in Z$ and let $S\subset N^0|_p$ be a 2--sphere with center
at zero and oriented in the standard way as $S^2\subset\R^3$. Then the
restriction of $K$ to $S$ has first Chern class equal to 2, and so is
non-trivial.
\end{lemma}

{\bf Proof of Lemma \ref{lem2.1}}\qua For simplicity it is enough to
consider the  case where the two positive eigenvalues of $\tau\dt\nabla\o$
are equal to 1 in as much as the Chern class of $K$ is unchanged by
continuous deformations of $\o$ near $Z$ which leave $\o^{-1}(0)$ unchanged.
With this understood, one can choose oriented local coordinates $(t,x,y,z)$
near the given point
$p$ so that $p$ corresponds to the origin, $Z$ is the set where $x=y=z=0$
and $dt$ is positive on $Z$ with respect to the given orientation. In these
coordinates,
\begin{equation}
\label{eq2.5}
\o=dt\wedge (xdx+ydy-2zdz)+x \ dy\wedge dz-y \ dx\wedge dz-2z \ dx
\wedge dy .
\end{equation}
(In these coordinates, the line $N_1$ corresponds to the $z$ axis.)

\bs
The strategy will be to identify a section of $K$ with nondegenerate
zeros  on $S$ and compute the Chern class by summing the degrees of these
zeros. There are three steps to this strategy. The first step identifies
$T^{1,0}\subset TX\ot\C$ and for this purpose it proves convenient to
introduce the functions
\begin{itemize}
\item $f  = 2^{-1}(x^2+y^2-2z^2)$,
\item $h  = (x^2+y^2)z$,
\item $g  = (x^2+y^2+4z^2)^{\frac 12}$.\eqnum
\end{itemize}
%\begin{equation}
%\label{eq2.6}
%\end{equation}
Also, introduce the standard polar coordinates $\rho =(x^2+y^2)^{\frac 12}$
and $\phi={\mbox{Arctan}}(y/x)$ for the $xy$--plane. Then $\o$ can be
rewritten as
\begin{equation}
\label{eq2.7}
\o=dt\wedge df+d\phi\wedge dh .
\end{equation}
Moreover, $\{dt,g^{-1}df,\rho \ d\phi,(g\rho)^{-1}dh\}$ form an oriented
orthonormal frame.

With the preceding understood, it follows that $T^{1,0}$ is the span of
$\{w_0\equiv dt+ig^{-1}df, \ w_1\equiv\rho \ d\phi+i(g\rho)^{-1}dh\}$ where
$\rho\neq 0$.

The second step in the proof produces a convenient section of $K$. In
particular, with $T^{1,0}$ identified as above, then $w_0\wedge w_1$ defines
a section of the canonical bundle with constant norm, at least where
$\rho\neq 0$. However, this section is singular at $\rho =0$. But the
section $\rho w_0\wedge w_1$ is nonsingular, vanishes on the sphere $S$ only
at the north and south poles and has nondegenerate zeros. (Since $w_0\wedge
w_1$ has constant norm it follows that the zeros of $\rho w_0\wedge w_1$ on
$S$ occur only where $\rho =0$ and are necessarily nondegenerate.)

The final step in the proof computes the degrees of the zeros of
$\rho w_0\wedge w_1$. For this purpose, it is convenient to first digress to
identify the fiber of $K$ where $\rho =0$. To start the digression, remark
that near $\rho =0$, one has $w_0=dt-i\e dz$ and $w_1=\rho d\phi+i\e d\rho$
to order $\rho$, where $\e$ is the sign of $z$. This implies that
$T^{1,0}\subset TX\ot\C$ at $\rho =0$ is spanned by the forms $w_0 = dt -
i\e dz$ and $-i \e e^{-i\e\phi} w_1 = dx - i\e dy$.  The latter gives $K
\subset \L^2(TX \otimes \C)$ where $\rho = 0$ as the span of $w_0 \wedge (dx
- i\e dy)$.

End the digression.  With respect to these trivializations, the section
$\rho w_0 \wedge w_1$ near $\rho = 0$ behaves to leading order in $\rho$ as
$\sim i\e (x+i\e y)$.  This last observation implies that the chosen section
$\rho w_0 \wedge w_1$ has Chern number equal to $2$ on $S$ as claimed.

\bs{\bf d) Pseudo-holomorphic subvarieties in X-Z}

A submanifold $C \subset X - Z$ is pseudo-holomorphic when $J$ maps $TC$ to
itself.  Note that such submanifolds have a canonical orientation as the
form $\o$ restricts to $TC$ as a nowhere vanishing $2$--form.

Here is an equivalent definition of a finite energy, pseudo-holomorphic
subvariety of $X-Z$:  The latter, $C$, is characterized by the following
conditions:
\begin{itemize}
\item $C$ is closed.
\item There is a countable set, $\L \subset C$, without accumulation points
and such that $C-\L$ is a pseudo-holomorphic submanifold.
\item $\int_{C-\L} \o < \i$.\eqnum
\end{itemize}
%\begin{equation}
%\label{eq2.8}
%\end{equation}

Note that a finite energy, pseudo-holomorphic subvariety $C$ naturally
defines a relative homology class, $[C] \in H_2(X,Z;\Z)$. As $H^2(X-Z;\Z)$
is the Poincare dual of $H_2(X,Z;\Z)$, there is no natural extension of the
intersection pairing to $H_2(X,Z)$.

Here are some examples of a finite energy, pseudo-holomorphic subvarieties:
First, let $M$ be a compact, oriented $3$--manifold with $b^1 > 0$.  Choose a
non-zero class in $H^1(X;\Z)$ and find a metric on $M$ for which the chosen
class is represented by a harmonic $1$--form with transversal zeros.  (A
generic metric will have this property.  See, eg \cite{Ho}.)  Let $\nu$
denote
the harmonic $1$--form.  Also, let $*$ denote the Hodge star for the metric
on $M$.  Now, take $X = S^1 \x M$ with metric that is the sum of that on $M$
with the metric on $S^1$ determined by a Euclidean coordinate $t \in
[0,2\pi]$.  Then $\o = dt \wedge \nu + *\nu$ is a harmonic, self-dual
$2$--form on $X$ where $Z = S^1 \x \{\nu^{-1}(0)\}$.  To see a
pseudo-holomorphic subvariety, introduce a flow line, $\g$, for the
vector field which is dual to $\nu$.  (Thus, $*\nu$ annihilates $T_{\g}$.)
Then $C = S^1 \x \g$ is a finite energy, pseudo-holomorphic submanifold if
$\g$ is either diffeomorphic to a circle or else is a path in $M$ connecting
a pair of zeros of $\nu$.  Note that when $\g$ is a path which connects a
pair of zeros to $\nu$, then the resulting $C$ will have intersection number
$1$ with any linking $2$--sphere of the corresponding pair of components of
$Z$.

Another example for this same $X$ has $C = S^1 \x \cup_i \g_i$, where
$\{\g_i\}$ is a finite set of flow lines for $\nu$ with each $\g_i$ being
either a circle or a path in $M$ connecting a pair of zeros of $\nu$.  Note
that if the flow lines in the set $\{\g_i\}$ precisely pair the
zeros of $\nu$, then the resulting pseudo-holomorphic variety $S^1 \x
\cup_i \g_i$ has intersection number $1$ with each linking $2$--sphere of $Z$
as required by Theorem \ref{th1.2}.

\bs{\bf e) The existence of pseudo-holomorphic subvarieties}

This subsection states a more detailed existence theorem for finite energy,
pseudo-holomorphic subvarieties.

The existence theorem for pseudo-holomorphic subvarieties (Theorem
\ref{th2.2}, below) uses solutions for a particular version of (\ref{eq2.2})
to construct the subvariety.  Here is the appropriate version:  Fix a
$\mbox{Spin}^{\C}$ structure for $X$ and a number $r \ge 1$.  Suppose that
$\o$ is a non-trivial, self-dual, harmonic $2$--form on $X$ and consider the
equations:
\begin{itemize}
\item $D_A\psi = 0$
\item $F_A^+ = rq(\psi)-i4^{-1}r\o$\eqnum
\end{itemize}
%\begin{equation}
%\label{eq2.9}
%\end{equation}
for a pair $(A,\psi)$ consisting of a connection $A$ on the line bundle $L =
\mbox{det}(S_+)$ and a section $\psi$ of $S_+$.  Note that these equations
constitute a version of (\ref{eq2.2}), as can be seen by replacing $\psi$
here with $\psi/\sqrt{r}$.

The precise statement of Theorem \ref{th2.2} requires the following three
part digression to explain some terminology.  Part~1 of the digression
introduces a numerical invariant for $\mbox{Spin}^{\C}$ structures.  For
this purpose, let $X$ be a compact, oriented, Riemannian $4$--manifold and
suppose that $\o$  a self-dual, harmonic $2$--form on $X$.  When $s$ is a
$\mbox{Spin}^{\C}$ structure for $X$, let $e_{\o}(s) \in \R$ denote the
evaluation on the fundamental class of $X$ of the cup product of $c_1(L)$
with the cohomology class of $\o$.

Part~2 of the digression introduces the scalar curvature $R_g$ for the
metric $g$, and also the metric's self-dual Weyl curvature, $W_g^+$.  These
tensors can be defined as follows:  Since $\L^2T^*X$ is the bundle
associated to the frame bundle via the adjoint representation on the Lie
algebra of $SO(4)$, the Riemann curvature tensor canonically defines an
endomorphism of $\L^2T^*X$.  Moreover, with respect to the splitting
$\L^2T^*X = \L_+ \oplus \L_-$, this tensor has a $2 \x 2$ block form and the
upper block gives an endomorphism of $\L_+$.  The trace of the latter is
$R_g/4$ and the traceless part is $W_g^+$.  (See eg \cite{AHS}.)  Also,
introduce the volume form, $\mbox{dvol}_g$, for the metric $g$.

Part~3 of the digression introduces the notion of an {\it irreducible
component} of a pseudo-holomorphic subvariety.  To appreciate the
definition, remember that there is a countable set $\L \subset C$ with no
accumulation points and with the property that $C - \L$ is a submanifold.
With this understood, an irreducible component of $C$ is the closure of a
component of $C - \L$.

End the digression.

\setcounter{theorem}{1}
\begin{theorem}
\label{th2.2}
Let $X$ be a compact, oriented, Riemannian $4$--manifold with $b^{2+} \ge 1$
and let $\o$ be a self-dual, harmonic $2$--form on $X$ which vanishes
transversely.  Fix a $\mbox{Spin}^{\C}$ structure $s$ and suppose that
there exists an unbounded sequence $\{r_n\} \subset [1,\i)$ with the
property that for each $n$, that $r = r_n$ version of (2.9) has a
solution, $(A_n,\psi_n)$.  Then, there exists a finite energy,
pseudo-holomorphic subvariety $C \subset X-Z$ with the following properties:
\begin{itemize}
\item Let $\{C_a\}$ denote the set of irreducible components of $C$.  Then
there exists a corresponding set of positive integers $\{m_a\}$ such that
$2\Sigma_am_a[C_a] \in H_2(X,Z;\Z)$ is Poincar\'e dual to the first Chern
class of the line bundle $L \otimes K|_{X-Z}$.  In particular, this implies
that $C$ has intersection number equal to $1$ with each linking $2$--sphere
of $Z$.
\item $\int_C \o \le \zeta e_{\o}(s) + \zeta \int_X |\o|(|R_g| +
|W_g^+|)\mbox{dvol}_g$.  Here, $\zeta$ is a universal constant.
\item For each $n$, let $\a_n \equiv (4i)^{-1}(cl^+(\o) + 2i)\psi_n$ and let
$\Sigma_n \equiv \a_n^{-1}(0)$.  Then
\begin{equation}
\label{eq2.10}
\lim_{n \to \i} \{\sup_{x \in C} \mbox{\rm\ dist}(x,\Sigma_n) + \sup_{x \in
\Sigma_n} \mbox{\rm\ dist}(x,C)\} = 0.
\end{equation}
\end{itemize}
\end{theorem}

Note that Theorem \ref{th1.2} follows directly from Theorem \ref{th2.2}
given the first point in (2.4).

There are also versions of Theorem 2.2 which holds when $\o$ does not
vanish transversely. Here is the simplest of these versions:

\begin{theorem}
\label{th2.3}
Let $X$ be a compact, oriented, Riemannian $4$--manifold with $b^{2+} \ge 1$
and let $\o$ be a self-dual, harmonic $2$--form on $X$. Fix a
${\mbox{Spin}}^{\C}$ structure $s$ with non-zero Seiberg--Witten invariant.
Then there exists a finite energy, pseudo-holomorphic subvariety
$C\subset X-Z$ with the following properties
\begin{itemize}
\item Let $\{C_a\}$ denote the set of irreducible components of $C$.
Then there exists a corresponding set of positive integers
$\{m_a\}$ such that $2\sum_a m_a[C_a]\in H_2(X,Z;\Z)$ is Poincar\'e dual to
the first Chern class of the line bundle $L\otimes K|_{X-Z}$.
\item $\int_C\o\leq\zeta e_{\o}(s)+\zeta\int_X |\o|(|R_g|+|W^+_g|)d
{\mbox{vol}}_g)$. Here, $\zeta$ is a universal constant.
\end{itemize}
\end{theorem}

Remark that Theorem \ref{th2.3} makes no assumptions about the
structure of the zero set of $\o$, but its assumption of a non-zero
Seiberg--Witten invariant is more restrictive than the assumption that
(2.9) has solutions for an unbounded set of $r$ values.
However, a version of Theorem \ref{th2.3} with the latter assumption
can be proved using the techniques in the subsequent sections if some
mild restrictions are assumed about the degree of degeneracy of the
zeros of $\o$. For example, the conclusions of Theorem \ref{th2.3}
hold if it is assumed that (2.9) has solutions for an
unbounded set of $r$ values, and if it is assumed that $\o^{-1}(0)$ is
non-degenerate except at some finite set of points. In any event,
these generalizations of Theorem \ref{th2.3} will not be presented
here.

\bs The remainder of this article is occupied with the proofs of
Theorem \ref{th2.2} and Theorem \ref{th2.3}. In this regard, the reader
should note
that
Theorem \ref{th2.3} is essentially a corollary of Theorem \ref{th2.2} and a
non-compact version of Gromov's compactness theorem \cite{Gr}. Meanwhile,
the proof of
Theorem \ref{th2.2} mimics as much as possible that of Theorem 1.3 in \cite{T1}
which asserts the equivalent result in the case where $\o$ has no zeros. In
particular, some familiarity with the arguments in sections 1--6 of
\cite{T1} and the revised reprint of the same article in \cite{T6} will
prove helpful. (The revisions of \cite{T1} in \cite{T6} correct some minor
errors in the original.) The final arguments for Theorem \ref{th2.2} are
given in
Section 7 below.

\bs{\bf f) The proof of Theorem \ref{th2.3}}

If $s$ has non-zero Seiberg--Witten invariant, then there is a $C^2$
neighborhood of the given metric
on $X$ such that for all $r$ sufficiently large and for each metric in this
neighborhood, the corresponding version of (2.9) has a solution. With
this understood, take a sequence of metrics $\{g_\nu\}_{\nu=1,2,\dots}$
with the following properties
\begin{itemize}
\item For each $g_{\nu}$, there is a self-dual, harmonic form $\o_{\nu}$
which vanishes transversely.
\item The sequence $\{g_{\nu}\}$ converges to the given metric $g$ in
the $C^{\i}$ topology as $\nu\to\i$.
\item The corresponding sequence $\{\o_{\nu}\}$ converges to $\o$.\eqnum
\end{itemize}
%\begin{equation}
%\label{eq2.12}
%\end{equation}
The existence of such a sequence can be proved as in \cite{Ho} or
\cite{LeB}.

Now invoke Theorem \ref{th2.2} for each metric $g_{\nu}$ and the
corresponding form
$\o_{\nu}$. Theorem \ref{th2.2} produces for each index $\nu$ a finite energy,
pseudo-holomorphic
variety $C_{\nu}$. Moreover, Theorem \ref{th2.2}  finds a uniform energy
bound for
each $C_{\nu}$.

Theorem \ref{th2.3}'s subvariety $C$ is now obtained as a limit of the sequence
$\{C_{\nu}\}$. The limit is found via a non-compact version of the Gromov
compactness theorem (in \cite{Gr}) for pseudo-holomorphic curves. The
particular
non-compact version is given by Proposition 3.8 in \cite{T2}. (The compact
case of Gromov's
compactness theorem is discussed in detail by numerous authors,
for example \cite{Pa},\cite{PW},\cite{MS} and \cite{Ye}.)

\section{Integral and pointwise bounds}
\label{sec3}

The geometric context for this section is as follows: Here $X$ is a compact,
oriented, Riemannian 4--manifold with a self-dual harmonic 2--form $\o$ which
vanishes
transversely. Also assume that $|\o|\leq\sqrt{2}$ everywhere. Now fix a
Spin${}^{\C}$ structure $s$ and for $r\geq 1$, introduce the perturbed
Seiberg--Witten equations as in (2.9).

The purpose of this section is to establish some basic properties of a
solution $(A,\psi)$ of this version of (2.9). For the most part,
these estimates are local versions of the estimates in (1.24) of \cite{T1}
and the arguments below are more or less modified versions of those from
section 2 of \cite{T1}.

Before turning to the details, please take note of the following
convention: The Greek letter $\z$ will be used to represent the ``generic"
constant in as much as its value may change each time it appears.
One should imagine a suppressed index on $\z$ which numbers its
appearances. Unless otherwise stated, the value of $\z$ is independent of
any points in question and, furthermore, $\z$ is always independent of the
parameter $r$ which appears in (\ref{eq2.10}). In general, $\z$ depends only
on the chosen Spin${}^{\C}$ structure and on the Riemannian metric.

A similar convention holds for the symbol $\z_{\d}$ when $\d$ is given as a
specified (minimum) distance to $Z$. That is, various inequalities below
will be proved under an assumption that the distance to $Z$ is greater than
some a priori value, $\d$. In these equations, $\z_{\d}$ denotes a generic
constant which depends only on the Spin${}^{\C}$ structure $s$, the
Riemannian metric, and the given $\d$. Furthermore, the precise value of
$\z_{\d}$ is allowed to change each time it appears and so the reader
should assume that $\z_{\d}$, like $\z$, is implicitly labeled by the order
of its appearance.

\bs{\bf a) Integral bounds for} $|\psi|^2$

The purpose of this subsection is to obtain pointwise estimates for the
components of the spinor $\psi$. The basis for these estimates is the
Bochner-Weitzenboch formula for $D_A^*D_A$ which, when applied to $\psi$,
reads
\setcounter{equation}{0}
\begin{equation}
\label{eq3.1}
\nabla_A^*\nabla_A\psi +4^{-1}R_g\psi+2^{-1}c_+(F^+_A)\psi = 0 ,
\end{equation}
where $R_g$ is the scalar curvature of the Riemannian metric. The strategy
below uses this last equation to generate first integral bounds for $\psi$
and then pointwise bounds.

The statement below of these integral bounds requires the reintroduction of
the notation of Theorem \ref{th2.2}. Here are the promised integral bounds:

\begin{lemma}
\label{lem3.1}
There is a universal constant $c$ with the following significance:
Let $s$ be a ${\mbox{Spin}}^{\C}$ structure on $X$. Now suppose that
$(A,\psi)$ solve (2.9) for the given ${\mbox{Spin}}^{\C}$ structure
and for some $r\geq 1$. Then
\begin{equation}
\label{eq3.2}
\begin{array}{rl}
\displaystyle{\int_X (2^{-\frac 12}|\o|-|\psi|^2)^2} &\leq\ 
\displaystyle{cr^{-1}(e_{\o}(s)+\int_X (|R_g|+r^{-1}|R_g|^2)
|\o|d{\mbox{vol}}_g) } \\
\displaystyle{} &  \\
\displaystyle{\int_X |\o| \ ||\psi|^2-2^{-\frac 12}|\o| |} &\leq \
\displaystyle{cr^{-1}(e_{\o}(s)}\\
&+\displaystyle{\int_X ((|R_g|+|W_g^+|) |\o|+r^{-1}|R_g|^2)
d{\mbox{vol}}_g).}
\end{array}
\end{equation}
\end{lemma}

{\bf Proof of Lemma \ref{lem3.1}}\qua The proof starts with the
observation that
$c_1(L)$ is represented in DeRham cohomology by $i\dt(2\pi)^{-1}F_A$
and thus $2\pi e_{\o}(s)=\int_X iF_A\wedge\o$. And, as $\o$ is self-dual,
the latter integral is equal to that of $iF_A\wedge\o$. Thus, line two of
(2.9) implies that
\begin{equation}
\label{eq3.3}
4\pi e_{\o}(s)\geq 2^{-\frac 12} r\int_X |\o|( 2^{-\frac 12}
|\o|-|\psi|^2) .
\end{equation}

To proceed, contract both sides of (\ref{eq3.1}) with $\psi$. The resulting
equation implies a differential inequality which can be written first as
\begin{equation}
\label{eq3.4}
2^{-1} d^*d|\psi|^2 +4^{-1}r|\psi|^2(|\psi|^2-
 2^{-\frac 12}|\o|)+4^{-1}R_g|\psi|^2\leq 0  .
\end{equation}
To put this last inequality in useful form, note first that rewriting
$|\psi|^2$ in its second and fourth appearances produces
\begin{equation}
\label{eq3.5}
\begin{array}{l}
\displaystyle{2^{-1} d^*d|\psi|^2 +(4\sqrt{2})^{-1} r(|\o|
(|\psi|^2-2^{-\frac 12}|\o|)+4^{-1}r(|\psi|^2-2^{-\frac 12}|\o|)^2}\\
\quad \displaystyle{+
(4\sqrt{2})^{-1} R_g|\o|+4^{-1}R_g (|\psi|^2-2^{-\frac 12}|\o|)\leq 0 .}
\end{array}
\end{equation}
Then an application of the triangle inequality to (\ref{eq3.5}) yields
\begin{equation}
\label{eq3.6}
\begin{array}{l}
\displaystyle{2^{-1} d^*d|\psi|^2 +(4\sqrt{2})^{-1} r|\o|
(|\psi|^2-2^{-\frac 12}|\o|)+8^{-1}r(|\psi|^2-2^{-\frac 12}|\o|)^2}\\
\qquad\qquad\qquad \displaystyle{\leq
-(4\sqrt{2})^{-1} R_g|\o|+8^{-1}r^{-1}R_g^2 .}
\end{array}
\end{equation}
Integrate this last inequality and compare with (\ref{eq3.3}) to obtain the
first line in (\ref{eq3.2}).

To obtain the
second line in (\ref{eq3.2}), note that the Weitzenboch formula for the
harmonic form $\o$ (see, eg~Appendix C in \cite{FU}) implies that
\begin{equation}
\label{eq3.7}
d^*d|\o|+ |\o|^{-1} |\nabla\o|^2\leq c|\cal W^+| \ |\o| ,
\end{equation}
where $\cal W^+$ is a universal curvature endomorphism of $\L_+$ which is
constructed from $R_g$ and $W^+_g$. In any event, with (\ref{eq3.7})
understood, introduce $u\equiv |\psi|^2-2^{-\frac 12}|\o|$. It now follows
from (\ref{eq3.7}) and (\ref{eq3.2}) (with an application of the triangle
inequality) that
\begin{equation}
\label{eq3.8}
2^{-1}d^*du+4^{-1}r|\o|u \leq c( |\o|(|R_g|+|W^+_g|)+|\o|^{-1}
|\nabla\o|^2 +r^{-1}|R_g|^2) .
\end{equation}

Now, integrate this last equation over the domain $\O\subset X$ where
$u\geq 0$ and then integrate by parts to find that
\begin{equation}
\label{eq3.9}
\int_X
|\o|u (|\psi|^2-2^{-\frac 12}|\o|)_+
\leq cr^{-1}\int_X (|\o|(|R_g|+|W^+_g|)+|\o|^{-1}
|\nabla\o|^2 +r^{-1}|R_g|^2) .
\end{equation}
Here $(|\psi|^2-2^{-\frac 12}|\o|)_+$ is the maximum of zero and
$|\psi|^2-2^{-\frac 12}|\o|$. (The boundary term which appears on the
left-hand side from integrating $d^*du$ is non-negative. This is easiest to
see when 0 is a regular value of $u$, for in this case the boundary
integral is minus that of the outward pointing normal derivative of $u$.
And, minus the latter derivative is non-negative as $u\geq 0$ inside $\O$
and $u\leq 0$ outside.)

With regard to (\ref{eq3.9}), note that $|\o|^{-1}$ is integrable across $Z$
since $|\o|$ near $Z$ is bounded from below by a multiple of the distance
to $Z$. Moreover, the integral over $X$ of $|\o|^{-1}|\nabla\o|^2$ can be
evaluated by integrating both sides of (\ref{eq3.7}). In particular, an
integration by parts eliminates the $d^*d|\o|$ integral, and one finds the
integral over $X$ of $|\o|^{-1}|\nabla\o|^2$ bounded by a universal
multiple of the integral over $X$ of $|\o|(|R_g|+|W_g|)$.

With this last point understood, the second line in (\ref{eq3.2}) follows
directly from (\ref{eq3.9}) and (\ref{eq3.3}).

\bs{\bf b) Pointwise bounds for} $|\psi|^2$

The purpose of this subsection is to derive pointwise bounds for
$|\psi|^2$. These bounds come from (\ref{eq3.4}) as well, but this time with
the help of the maximum principle. In particular, with the
maximum principle, (\ref{eq3.4}) immediately gives the bound
\begin{equation}
\label{eq3.10}
|\psi|^2\leq 1+r^{-1}\sup_X |R_g| .
\end{equation}
(Remember that $|\o|\leq\sqrt{2}$.) Here are some more refined bounds:

\begin{lemma}
\label{lem3.2} Let $\s(\dt)$ denote the distance function to $Z$.
There is a constant $\z$ which depends on the Riemannian metric
and is such that
\begin{equation}
\label{eq3.11}
\begin{array}{rll}
|\psi|^2 & \leq  & \displaystyle{2^{-1} +\z r^{-\frac 13}} \\
|\psi|^2 &\leq & 2^{-1}+\z r^{-1}\s^{-2} .
\end{array}
\end{equation}
\end{lemma}

The remainder of this section is occupied with the

\bs{\bf Proof of Lemma \ref{lem3.2}}\qua Start with the observation that the
right-hand side of (\ref{eq3.8}) is bounded by $\z\s^{-1}$, and thus
(\ref{eq3.8})  implies
\begin{equation}
\label{eq3.12}
2^{-1}d^*du+4^{-1} r|\o|u\leq \z\s^{-1} .
\end{equation}
To obtain the first bound in (\ref{eq3.11}), introduce a standard bump
function,
\begin{equation}
\label{eq3.13}
\chi\co [0,\i)\to [0,1] ,
\end{equation}
which is non-increasing, equals 1 on [0,1] and equals 0 on
$[2,\i)$. Given $R > 0$ and $r\geq 1$, promote $\chi$ to the function
$\chi_R\equiv\chi(r^{\frac 13}\s/R)$ on $X$. Then there is a constant
$\z_1$ such that the function $u'\equiv u+\z_1\chi_1\s$ obeys the
differential inequality
\begin{equation}
\label{eq3.14}
2^{-1}d^*du^{\prime}+4^{-1} r|\o|u^{\prime} \leq \z(1-\chi_{1/2})\s^{-1} .
\end{equation}
Then there is a constant $\z_2$ such that $u''=u+\z_1\chi_1\s-\z_2r^{-\frac
13}$ obeys
\begin{equation}
\label{eq3.15}
2^{-1}d^*du^{\prime\prime}+4^{-1} r|\o|u^{\prime\prime} \leq
-4r^{2/3} |\o|\z_2+   \z(1-\chi_{1/2})\s^{-1} \leq 0 .
\end{equation}
The previous equation and the maximum principle imply the first line in
(\ref{eq3.11}).

To obtain the second  line in
(\ref{eq3.11}), fix $c\geq 1$ and let $u'$ now denote $u-cr^{-1}\s^{-2}$. Then
(\ref{eq3.12}) implies that $u'$ obeys the differential
inequality
\begin{equation}
\label{eq3.16}
2^{-1}d^*du'+4^{-1} r|\o|u'\leq
-c4^{-1}\z^{-1}\s^{-1}+\z\s^{-1}+\z cr^{-1}\s^{-4} .
\end{equation}
(This is because $|\o|\geq\z^{-1}\s$ and because
$|\nabla^m\s|\leq\z_m\s^{1-m}$.) It follows from (\ref{eq3.16}) that if $c$
is taken large $(c\geq 8\z^2)$ and if $\s\geq 2z^{2/3}r^{-1/3}$, then
$2^{-1}d^*du'+4^{-1}r|\o|u'\leq 0$. Thus, for $c$ so constrained, the
function $u'$ cannot take a positive maximum where $\s\geq\z r^{-1/3}$.
Moreover the first line of (\ref{eq3.11}) insures that $u\leq \z r^{-1/3}$
where $\s=\z r^{-1/3}$, and so for $c\geq\z^3$, the function $u'$ cannot be
positive at all where $\s=\z r^{-1/3}$. This statement implies the second
line of (\ref{eq3.11}).

\bs{\bf c) Writing $\psi=(\a,\b)$ and estimates for} $|\b|^2$

To proceed from here, it proves convenient to introduce the components
$(\a,\b)$ of $\psi$ as follows:
\begin{equation}
\label{eq3.17}
\begin{array}{rll}
\a & \equiv & 2^{-1}(1+i(\sqrt{2} |\o|)^{-1}c_+(\o))\psi \\
\b & \equiv & 2^{-1}(1-i(\sqrt{2} |\o|)^{-1}c_+(\o))\psi
\end{array}
\end{equation}
The estimates in this subsection will show that $|\b|$ is uniformly small
away from $Z$. The following lemma summarizes

\begin{proposition}
\label{prop3.3}
Fix $\d>0$ and the ${\mbox{Spin}}^{\C}$ structure $s$. There are constants
$c_\d$, $\z_{\d}\geq 1$ which depend only on $s$, $\d$ and the Riemannian
metric and have the following significance:  Suppose that $(A,\psi)$ is a
solution to (2.9) as defined  by $s$ and with $r\geq \z_{\d}$. Then
the $\b$ component of $\psi$ obeys
\begin{equation}
\label{eq3.18}
|\b|^2\leq c_{\d}r^{-1}(2^{-1/2}|\o|-|\a|^2)+\z_{\d}r^{-2}
\end{equation}
at all points of $X$ with distance $\d$ or greater from $Z$.
\end{proposition}

The remainder of this subsection is occupied with the

\bs{\bf Proof of Proposition \ref{prop3.3}}\qua The estimates for the norm
of $\b$ are
obtained using the maximum principle with the projections
$2^{-1}(1+i(\sqrt{2}|\o|)^{-1}c_+(\o))$ of (\ref{eq3.1}). To start,
take the inner product of (\ref{eq3.1}) with the spinors $(\a,0)$ and
then $(0,\b)$ to obtain the following schematic equations:

\begin{equation}
\label{eq3.19}
\begin{array}{l}
\displaystyle{2^{-1}d^*d|\a|^2
+|\nabla_A\a|^2+4^{-1}r|\a|^2(|\a|^2-2^{-1/2}|\o|+|\b|^2)}
\qquad\qquad\qquad \\
\displaystyle{\qquad\qquad\qquad\qquad+{\cal R}_1(\a,\a)+{\cal R}_2(\a,\b)+{\cal
R}_3(\a,\nabla_A\b) = 0}
\\ \\
\displaystyle{2^{-1}d^*d|\b|^2 +|\nabla_A\b|^2
+4^{-1}r|\b|^2(|\b|^2-2^{-1/2}|\o|+|\a|^2)} \qquad \\
\displaystyle{\qquad\qquad\qquad\qquad +{\cal P}_1(\b,\b)+{\cal P}_2(\b,\a)+{\cal
P}_3(\b,\nabla_A\a) =0}
\end{array}
\end{equation}

Here, $\{{\cal R}_j\}$ and $\{{\cal P}_j\}$ are metric dependent endomorphisms.
Also, the covariant derivative (denoted above by $\nabla_A$) on the $\pm
\sqrt{2}i|\o|$ eigenbundles of the action of $c_+(\o)$ on $S_+$ are obtained by
projecting the spin covariant derivative on
$C^{\i}(S_+)$.

With regard to (\ref{eq3.19}), the associated connection for the derivative
$\nabla_A$ on sections of the $-\sqrt{2}i|\o|$ eigenbundles of $c_+(\o)$ is, after
an appropriate bundle identification, equal to half the difference between
$A$ and a certain canonical connection on the line bundle $K^{-1}$.
In particular, the associated curvature 2--form for this covariant
derivative is half the difference between $F_A$ and the
curvature of the canonical connection on  $K^{-1}$.

The canonical connection on  $K^{-1}$ is defined as in \cite{T4} and in
Section 1c of \cite{T1} as follows: There is a unique (up to isomorphism)
${\mbox{Spin}}^{\C}$ structure for $X-Z$ with the property that the
$-\sqrt{2}i|\o|$ eigensubbundle for the
$c_+(\o)$ action on the corresponding $S_+$ is a trivial complex line
bundle. For this ${\mbox{Spin}}^{\C}$ structure,the corresponding line
bundle $L$ is isomorphic to
$K^{-1}$. Thus there is a unique connection (up to gauge equivalence) on
$K^{-1}$ which makes the induced connection on the $-\sqrt{2}i|\o|$ subbundle
trivial. The latter connection is the canonical one. An alternate
definition of the covariant derivative of the canonical connection  on
$K^{-1}$ uses the natural identification of $K^{-1}$ (as an $\R^2$ bundle
over $X-Z$) as the orthogonal complement in $\L_+$ of the span of $\o$.
With this identification understood, the covariant derivative of the
canonical connection on a section of $K^{-1}$ is the orthogonal projection
onto $K^{-1}\ot T^*X$ of the Levi-Civita covariant derivative on
$C^{\i}(\L_+)$.

Likewise the associated curvature 2--form for the covariant derivative
$\nabla_A$ on sections of the $+\sqrt{2}i|\o|$ eigenbundle of $c_+(\o)$ is half the
sum of the curvatures of $A$ and the canonical connection on
$K^{-1}$. Moreover, after the appropriate bundle identifications, the
associated connection is, in a certain sense, half the sum of $A$ and the
canonical connection.

In any event, with (\ref{eq3.19}) understood, fix $\d >0$, but small so that it
is a regular value for the function $\s$. Then consider the second line of
(\ref{eq3.19}) where the distance to $Z$ is larger than $\d/2$. In particular,
where $\s\geq \d/2$ and when $r$ is large $(r\d >>1)$, the equation in the
second line of (\ref{eq3.19}) for $|\b|^2$ yields the inequality
\begin{equation}
\label{eq3.20}
2^{-1}d^*d |\b|^2+|\nabla_A\b|^2+(4\sqrt{2})^{-1}r|\o|
|\b|^2\leq \z r^{-1}\d^{-1}(|\a|^2+|\nabla_A\a|^2) .
\end{equation}
To use this last equation, write $w=(2^{-1/2}|\o|-|\a|^2)$ and
then observe that (\ref{eq3.7}) and the first line of (\ref{eq3.19}) yield
(where $\s\geq\d/2$) the inequality
\begin{equation}
\label{eq3.21}\begin{split}
2^{-1}d^*d (-w)+|\nabla_A\a|^2+&\,(4\sqrt{2})^{-1}r|\o|(-w)+4^{-1} rw^2\\
\leq\,&\z (|\a|^2+|\nabla_A\b|^2+|\b|^2+\d^{-1}).\end{split}
\end{equation}

Add a large (say $c_{\d}\geq 1$) multiple of $r^{-1}$ times this last
equation to (\ref{eq3.7}) to obtain the following inequality at points where
$\s\geq \d/2$
\begin{equation}
\label{eq3.22}
2^{-1}d^*d (|\b|^2-c_{\d}r^{-1}w)+(4\sqrt{2})^{-1}r|\o|
(|\b|^2-c_{\d}r^{-1}w)
\leq \z_{\d} r^{-1} .
\end{equation}
(Remember the convention that $\z_{\d}$ is a constant which depends only on
the ${\mbox{Spin}}^{\C}$ structure $\d$ and on the Riemannian metric.
Furthermore, the precise value
of $\z_{\d}$ can change from  appearance to appearance.)
In particular, there is an $r$--independent constant $\z_{\d}\geq 1$ which
guarantees (\ref{eq3.22}) when $c_{\d}\in (\z_{\d},\z_{\d}^{-1}r)$.
With the preceding understood, (\ref{eq3.22}) implies that there exists a
constant $\z$ which is independent of both $r$ and $\d$, and is such that
the function $u\equiv|\b|^2-c_{\d}r^{-1}w-4\sqrt{2}r^{-2}\z\z_{\d}$ obeys
\begin{equation}
\label{eq3.23}
2^{-1}d^*du +(4\sqrt{2})^{-1}r|\o| u \leq 0
\end{equation}
at all points where $\s\geq \d/2$ when $r$ is large.

Now, hold (\ref{eq3.23}) for the moment and consider that there is a unique,
continuous function $v$ which equals 1 where $\s\leq\d/2$ and satisfies
$2^{-1}d^*dv+(4\sqrt{2})^{-1}r|\o|v=0$ where $\s\geq\d/2$. Furthermore, $v$ is
positive and, as $|\o|\geq\z^{-1}\d$, this $v$ obeys
\begin{equation}
\label{eq3.24}
v \leq \z\exp(-\sqrt{r} (\s-\d)/\z_{\d}) .
\end{equation}
To see this last bound, note first that the maximum principle implies that
$v$ is no greater than the function $v'$ which is 1 where $\s\leq\d/2$ and
obeys $2^{-1}d^*dv'+(4\sqrt{2})^{-1}r\z^{-1}\d v'=0$ where $\s\geq\d/2$.
The value of $v'$ at $x$ can be written as an integral of the derivative of
a Green's kernel $G(x,\dt)$ over the boundary set where $\s=\d/2$.
However, when $m\geq 0$ is a constant, the Green's kernel $G(x,y)$ for the
operator $d^*d+m^2$ obeys the bound
\begin{equation}
\label{eq3.25}
\z^{-1}{\mbox{ dist}}(x,y)^{-2}e^{-\z m{\mbox{\footnotesize\ dist}}(x,y)}\leq
|G(x,y)| \leq \z {\mbox{ dist}}(x,y)^{-2}
e^{-m{\mbox{\footnotesize\ dist}}(x,y)/\z} .
\end{equation}
This last bound yields the bound for $v'$ by the right-hand side of
(\ref{eq3.24}).

Now, as Lemma \ref{lem3.2} insures that $u\leq\z\d$ where
$\s=\d/2$, it follows from the maximum principle that $u-\z\d v$ is
non-positive where $\s\geq\d/2$. This last estimate implies that
$|\b|^2\leq c_{\d}r^{-1} w+\z_{\d}r^{-2}$ at points where $\s\geq\d$
which is the statement of Proposition \ref{prop3.3}.

\bs{\bf d) Bounds for the curvature}

This subsection modifies the arguments in section 2d of \cite{T1} to
establish bounds for the curvature $F_A$ of the connection part of
$(A,\psi)$. In this regard, the estimates for the self-dual part, $F^+_A$,
come directly from (\ref{eq3.10}) and the defining equation for $F^+_A$ in the
second line of (2.9). In particular, the second line of (2.9)
implies that
\begin{equation}
\label{eq3.26}
|F^+_A|=r(2\sqrt{2})^{-1} ((2^{-1/2}|\o|-|\a|^2)^2 +
2|\b|^2(2^{-1/2}|\o|+|\a|^2)+|\b|^4)^{1/2} .
\end{equation}
Moreover, with $\d>0$ specified and with $r$ large, Proposition \ref{prop3.3}
implies the more useful bound
$|F^+_A|=r(2\sqrt{2})^{-1} ((2^{-1/2}|\o|-|\a|^2)+r^{-1}
\z'_{\d})^2+r^{-2}\z^2_{\d})^{1/2}$ at points where the distance,
$\s$, to $Z$ is greater than $\d$. Thus the triangle inequality gives
\begin{equation}
\label{eq3.27}
|F^+_A|\leq r(2\sqrt{2})^{-1} (2^{-1/2}|\o|-|\a|^2)+\z_{\d}
\end{equation}
at each point where $\s\geq\d$ and when $r\geq\z_{\d}$.

With the preceding understood, consider next the case of
$|F^-_A|$. The following proposition summarizes:

\begin{proposition}
\label{prop3.4}
Fix a ${\mbox{Spin}}^{\C}$ structure $s$, and fix $\d >0$. Then there are
constants
$\z_{\d}$, $\z'_{\d}\geq 1$ with the following significance:
Let $(A,\psi)$ be a solution to the $s$ version of (2.9)
where $r\geq\z_{\d}$. Then at points in $X$ with distance $\d$ or more from
$Z$,
\begin{equation}
\label{eq3.28}
|F^-_A|\leq r(2\sqrt{2})^{-1} (1+\z_d r^{-1/2})
 (2^{-1/2}|\o|-|\a|^2)+\z'_{\d}
\end{equation}
\end{proposition}

The remainder of this subsection is occupied with the

\bs{\bf Proof of Proposition \ref{prop3.4}}\qua The proof is divided into
seven steps.

{\bf Step 1}\qua This first step states and then proves a bound
on the $L^2$--norm of $F^-_A$ over the whole of $X$. The following lemma
gives this bound plus a second bound for $|\nabla_A\psi|$ and $F^+_A$ which
will be exploited in a subsequent step.

\begin{lemma}
\label{lem3.5}
There is a constant $\z$ which depends only on the metric and on
the ${\mbox{Spin}}^{\C}$ structure and which has the following significance:
Let $r\geq 1$ be given
and let $(A,\psi)$ be a solution to (2.9) using the
${\mbox{Spin}}^{\C}$ structure and $r$. Then
\begin{equation}
\label{eq3.29}
\begin{array}{l}
\displaystyle{\int_X |F_A|^2 \leq \z r} \\
\displaystyle{} \\
\displaystyle{\int_X (1+{\mbox{\rm\ dist}}(x,\dt)^{-2})(|\nabla_A\psi|^2+r^{-1}
|F_A^+|^2 ) \leq \z \quad {\mbox{for any point}} \ \ x\in X.}
\end{array}
\end{equation}
\end{lemma}

{\bf Proof of Lemma \ref{lem3.5}}\qua Take the inner product of both
sides of
(\ref{eq3.1}) with $\psi$ to obtain the equation
\begin{equation}
\label{eq3.30}
2^{-1}d^*d |\psi|^2+|\nabla_A\psi|^2+8^{-1}r^{-1} |F_A^+|^2\leq\z
|R_g|^2|\psi|^2+ir^{-1} \<F^+_A,\o\> .
\end{equation}
Here $\< \ , \ \>$ denotes the metric inner product on $\L^2T^*X$.
Integration of this last equation over $X$ yields a uniform bound on the
$L^2$ norms of $\nabla_A\psi$ and $F^+_A$. The first inequality in
(\ref{eq3.29}) then follows from the fact that the difference between the $L^2$
norms of $F^+_A$ and $F^-_A$ is equal to a universal multiple of the
evaluation of $c_1(L)\cup c_1(L)$ on the fundamental class of $X$.

To obtain the second inequality in (\ref{eq3.29}), introduce the Green's
function for the operator $d^*d+1$ with pole at $x$. Multiply both sides of
(\ref{eq3.30}) with $G(x,\dt)$ and integrate the result over $X$.
Then integrate by parts and use (\ref{eq3.25}) to obtain
\begin{equation}
\label{eq3.31}\begin{split}
2^{-1}|\psi|^2(x)+\z^{-1}&\int_X {\mbox{ dist}}(x,\dt)^{-2}
(|\nabla_A\psi|^2+r^{-1} |F_A^+|^2)\\&\leq\z
(r^{-1} +\int_X (|\nabla_A\psi|^2+r^{-1} |F_A^+|^2)\leq\z .
\end{split}
\end{equation}
This last equation implies the second line in (\ref{eq3.29}).

\medskip
{\bf Step 2}\qua This second step derives a differential
inequality for $|F^-_A|$. This derivation starts as in section 2d of
\cite{T1}. In particular, (2.14--2.15) in \cite{T1} hold in the present
case for $\mu=-2^{-1}iF^-_A$. (The factor of 2 here comes about because the
equations in \cite{T1} refers not to the connection $A$ on $L$, but to a
connection on a line bundle whose square is the tensor product of the
canonical bundle with $L$.) Now argue as in the derivation of \cite{T1}'s
(2.19), to find that $s\equiv |F^-_A|$ obeys the differential
inequality
\begin{equation}
\label{eq3.32}
\begin{array}{l}
2^{-1} d^*ds+4^{-1} r|\a|^2s\leq
\z s+(2\sqrt{2})^{-1} r(|\nabla_A\a|^2+|\nabla_A\b|^2) \\
\qquad \qquad \qquad \qquad \qquad + \z r(|\a|^2+|\b|^2+|\a| \ |\nabla_A\b| + |\b| \ |\nabla_A\a|) .
\end{array}
\end{equation}
%\begin{eqnarray*}
%\lefteqn{2^{-1} d^*ds+4^{-1} r|\a|^2s\leq
%\z s+(2\sqrt{2})^{-1} r(|\nabla_A\a|^2+|\nabla_A\b|^2) }\\
%& & + \z r(|\a|^2+|\b|^2+|\a| \ |\nabla_A\b| + |\b| \ |\nabla_A\a|) .
%\end{eqnarray*}
Note that this last equation holds everywhere on $X$. (Since $s$ is not
necessarily $C^2$ where $s=0$, one should technically interpret
(\ref{eq3.32}) as an inequality between distributions on the space of positive
functions. However, this and similar technicalities below have no bearing
on the subsequent arguments.  Readers who are uncomfortable with this
assertion can replace $s$ in (\ref{eq3.32}) and below by
$(|F_A^-|^2+1)^{1/2}$ without affecting the arguments.)

\medskip
{\bf Step 3}\qua This  step uses (\ref{eq3.32}) to bound $s$ by
$\z r$ everywhere on $X$ when $r\geq\z$. To obtain such a bound, multiply
both sides of (\ref{eq3.32}) by the Green's function $G(x,\dt)$ for the
operator $d^*d+1$. Integrate the resulting inequality over $X$ and
integrate by parts to obtain
\begin{equation}
\label{eq3.33}
s(x)+r\int_X |\a|^2 s\dt{\mbox{ dist}}(x,\dt)^{-2}\leq
\z\int_X s {\mbox{ dist}}(x,\dt)^{-2}+\z r .
\end{equation}
Here (\ref{eq3.25}) has been used. Also, the first line of (\ref{eq3.29})
has been
invoked to bound the integral over $X$ of $s$ by $r^{1/2}$. In addition,
the second line of (\ref{eq3.29}) has been invoked to bound the product of
$G(x,\dt)$ with $|\nabla_A\psi|^2$.

To make further progress, fix $R >0$ and break the integral on the right
side of (\ref{eq3.33}) into the part where dist$(x,\dt)\geq R$ and the
complementary region. With this done, (\ref{eq3.33}) is seen to imply that
\begin{equation}
\label{eq3.34}
\sup_X s+r\int_X |\a|^2 s\dt
 {\mbox{ dist}}(x,\dt)^{-2}\leq\z R^{-2}r^{1/2} +\z R^2\sup_X s+\z r .
\end{equation}
With $R=2^{-1} \z^{-1/2}$, this last line gives the claimed bound
\begin{equation}
\label{eq3.35}
\sup_X |F^-_A| \leq \z r .
\end{equation}

\medskip
{\bf Step 4}\qua This  step uses (\ref{eq3.32}) to derive a simpler
differential inequality. To begin, reintroduce $w\equiv
2^{-1/2}|\o|-|\a|^2$. It then follows from (\ref{eq3.20}) and
(\ref{eq3.21}) that
there are constants $\kappa_1$, $\kappa_2$ and $\kappa_3$ which depend only on the
metric (not on $(A,\psi)$ nor $r$) and which have the following
significance: Let $q_0\equiv (2\sqrt{2})^{-1} r(1+\kappa_1/r)w-\kappa_2 r|\b|^2+\kappa_3$
and then the function $(s-q_0)$ obeys $2^{-1}d^*d(s-q_0)+4^{-1} r|\a|^2(s-q_0)
\leq \z(s+r\s^{-1})$. Now introduce the function $(s-q_0)_+\equiv\max
((s-q_0),1)$. The latter function obeys the same differential
inequality as does $(s-q_0)$, namely
\begin{equation}
\label{eq3.36}
2^{-1} d^*d(s-q_0)_+ +4^{-1}r|\a|^2(s-q_0)_+\leq \z(s+r\s^{-1}) .
\end{equation}
(To verify (\ref{eq3.36}), write $(s-q_0)_+=2^{-1}((s-q_0)+|s-q_0|)$. Also,
since
(\ref{eq3.36}) involves two derivatives of the Lipschitz function
$|s-q_0|$, this last equation should be interpreted as an inequality
between distributions on the space of positive functions. As before, such
technicalities play no essential role in the subsequent arguments.)

\medskip
{\bf Step 5}\qua Now, fix $\d >0$ but small  enough to be a
regular value of $\s$. Let $\chi$ denote the bump function in (\ref{eq3.13})
and let $\chi^\d\equiv\chi(\s(\dt)/\d)$. Agree to let
$q_1\equiv (1-\chi^\d)(s-q_0)_+$. Note that the task now is to bound $q_1$
from above.

To begin, multiply both sides of (\ref{eq3.36}) by $(1-\chi^\d)$ to obtain
\begin{equation}
\label{eq3.37}
2^{-1} d^*dq_1 +4^{-1}r|\a|^2 q_1\leq \z_\d
r+\< d\chi^\d,d(s-q_0)_+\>
\end{equation}
at points where $\s\geq\d$. Here $\< \ , \ \>$ denotes the metric inner
product. Now observe that the maximum principle insures that
$q_1\leq q_2+q_3$, where $q_2$ solves the equation
\begin{equation}
\label{eq3.38}
2^{-1} d^*dq_2 +(4\sqrt{2})^{-1}r|\o|q_2= \z_\d
r+\< d\chi^\d,d(s-q_0)_+\> ,
\end{equation}
where $\s\geq\d$ and $q_2=0$ where $\s=\d$. In the mean time, $q_3$ obeys
\begin{equation}
\label{eq3.39}
2^{-1} d^*dq_3 +4^{-1}r|\a|^2 q_3= \z_\d
4^{-1}r|w| \ |q_2|
\end{equation}
where $\s\geq\d$ and vanishes where $\s=\d$. Here,
$w\equiv 2^{-1/2}|w\-|\a|^2$.

Bounds for $|q_2|$ can be found with help of the Green's function,
$G(\dt,\dt)$ for the operator $2^{-1}d^*d+(4\sqrt{2})^{-1}r|\o|$ with Dirichlet
boundary conditions on the surface $\s=\d$. In particular, since
$|\o|\geq\z^{-1}\d$ standard estimates bound $|\nabla^kG(x,y)|$ for $k=0,1$ at
points $x\neq y$ by $\z_\d|x-y|^{-2-k}\exp(-\sqrt{r}|x-y|/\z_\d)$.
(Since $|\o|\geq\zeta^{-1}\d$ where the distance to $Z$ is greater than $\d$,
these standard estimates involve little more than (\ref{eq3.25}) and the
maximum principle.)

Let $\D_\d$ denote the supremum of $(s-q_0)_+$ where $\s\geq\d$.
(Note that (\ref{eq3.25}) asserts that $\D_\d\leq\z_\d r$ in any event.)
Then the estimates just given for the Green's function imply that
\begin{equation}
\label{eq3.40}
 |q_2|\leq \z_\d (1+r^{-1/2} \D_\d \exp [-\sqrt{r} (\s-\d)/\z_\d ])
\end{equation}
at points where $\s\geq\d$ and when $r\geq\z_\d$.

\medskip
{\bf Step 6}\qua The purpose of this step is to obtain a bound
for the supremum norm of $|q_3|$. Such a bound is a part of the assertions
of

\begin{lemma}
\label{lem3.6}
Given $\d >0$, there is a constant $\z_{\d}\geq 1$ which is independent of
$(A,\psi)$ and $r$ and is such that if $r\geq\z_\d$, then the following is
true:
\begin{equation}
\label{eq3.41}
q_3\leq\z_\d(r^{1/2}+r^{-1/6}\D_\d) .
\end{equation}
\end{lemma}

This claim will be proved momentarily. Note however that (\ref{eq3.41}) leads
to a refinement of the bound $\D_\d\leq\z_\d r$ which came from (\ref{eq3.35}):
\begin{equation}
\label{eq3.42}
\D_\d\leq\z_\d r^{1/2}
\end{equation}
when $r\geq\z_\d$. To obtain this refinement, remark that according to
(\ref{eq3.40}) and (\ref{eq3.41}),
\begin{equation}
\label{eq3.43}
s-q_0\leq\z_\d(r^{1/2}+r^{-1/6}\D_\d)
\end{equation}
at points where $\s\geq 2\d$ when $r\geq\z_\d$. This last estimate implies
that $\D_{2\d}\leq\z_\d(r^{1/2}+r^{-1/6}\D_\d)$. The latter inequality,
iterated thrice, reads $\D_{8\d}\leq z_{8\d}r^{1/2}(1+r^{-1}\D_\d)$.
Now plug in the bound of $\D_\d$ by $\z_\d r$ from (\ref{eq3.35}) to conclude
that $\D_{8\d}\leq z_{8\d}r^{1/2}$. Replacing $8\d$ by $\d$ gives
(\ref{eq3.42}).

\bs{\bf Proof of Lemma \ref{lem3.6}}\qua Since $q_3\geq 0$, this function is no
greater than the solution, $u$, to the equation $2^{-1} d^*du=\z_\d
4^{-1} r|w|(1+r^{-1/2}\D_\d\exp[-\sqrt{r}(\s-\d)/\z_\d])$ where $\s\geq\d$
with Dirichlet boundary conditions where $\s=\d$. This function $u$ can be
bounded using the Green's function for the Laplacian. In particular
(\ref{eq3.25}) gives
\begin{equation}
\label{eq3.44}\begin{split}
u(x)&\leq\z_\d \int_{\s\geq\d} {\mbox{dist}}(x,\dt)^{-2}
r|w|\\& +\z_\d\D_\d  \int_{\s\geq\d} {\mbox{dist}}(x,\dt)^{-2}
r^{1/2} |w|\exp [-\sqrt{r}(\s-\d)/\z_\d]).\end{split}
\end{equation}
Consider the two integrals above separately. To bound the first integral,
fix $d>0$ but small and break the region of integration into the part where
dist$(x,\dt)\geq d$, and the complementary region. The integral over the
first region is no greater than that of $d^{-2}r|w|$ over the region where
$\s\geq\d$. Since the integral of $r|w|$ is uniformly bounded (Lemma
\ref{lem3.1})
by some $\z_\d$, this first part of the first integral in (\ref{eq3.44}) is no
greater than $\z_\d d^{-2}$.  Meanwhile, since $|w|\leq\z_\d$, the
dist$(x,\dt)\leq d$ part of the first integral above is no greater than
$\z_\d rd^2$. Thus, taking $d=r^{-1/4}$ bounds the first integral in
(\ref{eq3.44}) by $\z_\d r^{1/2}$.

Now consider the second integral, and again consider the contributions from
the region where dist$(x,\dt)\geq d$  and the complementary region. The
contribution from the first region is no greater than
$\z_\d\D_\d r^{-1/2}d^{-2}$ since $r|w|$ has a uniform bound on its
integral. Meanwhile the region where dist$(x,\dt)\geq d$ contributes no
more than $\z_\d\D_\d d$. Thus taking $d=r^{-1/6}$ bounds the second
integral in (\ref{eq3.44}) by $\z_\d\D_\d r^{-1/6}$.

\medskip
{\bf Step 7}\qua This step completes the proof of Proposition
\ref{prop3.4} with the help of a pointwise bound on $q_3$. To continue the
argument,
mimic the discussion surrounding $(2.27)$ and $(2.28)$ of \cite{T1}
to find constants $\z_\d$ and $c_\d$ such that the function
\begin{equation}
\label{eq3.45}
v_1\equiv w-c_\d |\b|^2 +\z_\d/r
\end{equation}
obeys the following properties at points where $\s\geq\d$ and
when $r\geq\z_\d$.
\begin{equation}
\label{eq3.46}
\begin{array}{rll}
v_1 &\geq & 2^{-1} \z_\d/r .\\
v_1 &\geq & w\\
2^{-1}d^*d v_1+4^{-1} r|\a|^2 v_1 &\geq & 0 .
\end{array}
\end{equation}

With $v_1$ understood, note that (\ref{eq3.41}) and (\ref{eq3.42}) and the
second
and third lines of (\ref{eq3.47}) imply that there exists $\z_\d\geq 1$
such that $q_4\equiv q_3-\z_\d r^{1/2}v_1$ obeys
\begin{equation}
\label{eq3.47}
2^{-1} d^*d q_4+4^{-1} r|\a|^2 q_4\leq 4^{-1} r|w| \ |q_2| \
\end{equation}
where $\s\geq\d$ and $q_4\leq 0$ where $|\a|^2\leq (2\sqrt{2})^{-1}|\o|$
and where $\s=\d$. This last point implies (via the maximum  principle)
that $q_3-\z_\d r^{1/2} v_1$ is no greater than the solution $v$ to the
differential equation $2^{-1}d^*dv+8^{-1}rv=4^{-1} r\z_\d \sup(|q_2|)$
where $\s\geq\d$ with boundary condition $v\geq 0$ where $\s=\d$.
In particular, it follows from (\ref{eq3.40}) and (\ref{eq3.42}) that
$v\leq\z_\d$ and thus
\begin{equation}
\label{eq3.48}
q_3-z_\d r^{1/2} v_1\leq \z_\d .
\end{equation}
This last bound with (\ref{eq3.40}), (\ref{eq3.42}) and (\ref{eq3.45})
complete the
arguments for Proposition \ref{prop3.4}.

\bs{\bf e) Bounds for $\nabla_A\a$ and} $\nabla_A\b$

This subsection modifies the arguments in section 2e of \cite{T1} to obtain
pointwise bounds on the covariant derivatives of $\a$ and $\b$. The
following proposition summarizes:

\begin{proposition}
\label{prop3.7}
Fix a ${\mbox{Spin}}^{\C}$ structure for $X$. Given $\d >0$, there are
constants $\z_\d$ and $\z'_\d$
with the following significance:  Let $r\geq\z_\d$ and let $(A,\psi)$ be a
solution to the $r$ version of (2.9) using the given
${\mbox{Spin}}^{\C}$ structure.
Then at points where the distance to $Z$ is larger than $\d$, one has
\begin{equation}
\label{eq3.49}
|\nabla_A\a|^2+r|\nabla_A\b|^2\leq\z_\d r(2^{-1/2}|\o|-|\a|^2)+\z'_d .
\end{equation}
\end{proposition}

The remainder of this section is occupied with the

\bs{\bf Proof of Proposition \ref{prop3.7}}\qua The arguments for
Proposition \ref{prop3.7} are slight modifications of those for
Proposition $2.8$ in
\cite{T1}. In any event, there are three steps to the proof.

\medskip
{\bf Step 1}\qua
For applications in a subsequent section, it proves convenient to introduce
$\underline{\a}\equiv 2^{1/2} |\o|^{-1/2}\a$. Note that an $r$--independent
bound for $|\nabla_A\underline{\a}|$ where the distance to $Z$ is greater
than $\d$ gives an $r$--independent
bound for $|\nabla_A\a|$.

The manipulations that follow assume that the distance to $Z$ is greater
than $\d$ and that $r\geq\z_\d$ so that Proposition \ref{prop3.3} and
Proposition \ref{prop3.4} can be invoked.

To begin, note that (\ref{eq3.1}) implies an equation for $\underline{\a}$
which has the following schematic form:
\begin{equation}
\label{eq3.50}
\nabla_A*\nabla_A\underline{\a}+4^{-1} r\kappa^2
\underline{\a}(|\underline{\a}|^2-1)+4^{-1} r\underline{\a} |\b|^2+
{\cal R}(\underline{\a},\b,\nabla_A\underline{\a},\nabla_A\b)=0 ,
\end{equation}
where $\cal R$ is multilinear in its four entries and satisfies
$|{\cal R}|+|\nabla{\cal R}|+|\nabla^2{\cal R}|\leq\z_\d$. Differentiate this
last equation and commute derivatives where appropriate to obtain
\begin{equation}
\label{eq3.51}\begin{split}
\nabla_A*\nabla_A(\nabla_A\underline{\a})&+4^{-1} r\kappa^2
\nabla_A\underline{\a} +{\cal Q}_1(
\underline{\a},\b)\\&+
{\cal Q}_2(\nabla_A^2\underline{\a},\nabla_A^2\b)+
{\cal T}_1\nabla_A\a +{\cal T}_2\nabla_A\b =0.\end{split}
\end{equation}
Here $\{{\cal Q}_j\}_{j=1,2}$ are bilinear in their entries.
Moreover, $|{\cal Q}_j|+|\nabla{\cal Q}_j|\leq\z_\d$ for $j=1$ and 2.
Meanwhile, $|{\cal T}_1|\leq\z_\d(1+r|w|)$ and $|{\cal T}_2|\leq\z_\d
(1+r|w|)^{1/2}$. (The latter bounds use Proposition \ref{prop3.3} and
Proposition
\ref{prop3.4}.

Next note that there is a similar equation for $\nabla_A\b$.
\begin{equation}
\label{eq3.52}\begin{split}
\nabla_A*\nabla_A(\nabla_A\b)&+4^{-1} r\kappa^2
\nabla_A\b +{\cal Q}'_1(
\underline{\a},\b)\\&+
{\cal Q}'_2(\nabla_A^2\underline{\a},\nabla_A^2\b)+
{\cal T}'_1\nabla_A\b +{\cal T}'_2\nabla_A\a =0 .\end{split}
\end{equation}
Here $\{{\cal Q}'_j\}$ and $\{{\cal T}_j\}$ obey the same bounds as their
namesakes in (\ref{eq3.51}).

\medskip
{\bf Step 2}\qua Take the inner product of (\ref{eq3.51}) with
$\nabla_A\underline{\a}$ and that of (\ref{eq3.52}) with $\nabla_A\b$.
Judicious use of the triangle inequality yields
\begin{equation}
\label{eq3.53}
\begin{array}{l}
\displaystyle{2^{-1}d^*d |\nabla_A\underline{\a}|^2 +
|\nabla_A\nabla_A\underline{\a}|^2 +(4\sqrt{2})^{-1} r|\o| |\nabla_A
\underline{\a}|^2} \\
\qquad \qquad \qquad \qquad \quad \displaystyle{\leq \z_\d (r^{-1} +(1+r|w|) |\nabla_A\a|^2
+|\nabla_A\b|^2) +r^{-1} |\nabla_A\nabla_A\b|^2} \\
\displaystyle{} \\
\displaystyle{2^{-1}d^*d |\nabla_A\b|^2 + |\nabla_A\nabla_A\b|^2
+(4\sqrt{2})^{-1} r|\o| |\nabla_A \b|^2} \\
\qquad \qquad \qquad \qquad \quad \displaystyle{\leq \z_\d (r^{-1} +(1+r|w|) |\nabla_A\b|^2
+|\nabla_A\a|^2)
+r^{-1} |\nabla_A\nabla_A\a|^2}
\end{array}
\end{equation}
Here $r\geq\z_\d$ is assumed so that Lemma \ref{lem3.2} and Proposition
\ref{prop3.3}
can be invoked.

Now introduce $y\equiv |\nabla_A\underline{\a}|^2+r|\nabla_A\b|^2$.
By virtue of (\ref{eq3.53}), the latter obeys
\begin{equation}
\label{eq3.54}
2^{-1}d^*d y+(4\sqrt{2})^{-1} r|\o| y\leq \z_\d
(1+r(|\nabla_A\a|^2 +r|\nabla_A\b|^2)) .
\end{equation}

\medskip
{\bf Step 3}\qua Reintroduce the function
$w\equiv 2^{-1/2} |\o|-|\underline{\a}|^2$. It then follows from (\ref{eq3.20})
and (\ref{eq3.21}) that there are constants $\kappa_{\d,1}$,  $\kappa_{\d,2}$ and
$\kappa_{\d,3}$
which depend only on $\d$ and are such that $y'\equiv y-\kappa_{\d,1}
rw+\kappa_{\d,2}r^2|\b|^2-\kappa_{\d,3}$ obeys
\begin{equation}
\label{eq3.55}
2^{-1}d^*d y'+(4\sqrt{2})^{-1} r|\o| y'\leq 0
\end{equation}
where the distance to $Z$ is greater than $\d$. Now introduce
$y'_+\equiv\max(y',0)$ and note that (\ref{eq3.51}) is still true with
$y'_+$ replacing $y'$, at least as a distribution on the space of positive
functions with support where the distance to $Z$ is greater than $\d$.

With the preceding understood, take the function $\chi$ from
(\ref{eq3.13}) and set $\chi^\d\equiv\chi(\s(\dt)/\d)$, where $\s$ denotes the
distance function to $Z$. Let $G(\dt,\dt)$ denote the Green's function for
the operator $d^*d+(4\sqrt{2})^{-1} r\z|\o|\d$, and with $x\in X$ obeying
$\s(x)\geq 2\d$. Here $\z$ is chosen so that $|\o|\geq\z\d$ at all points
where $\s\geq\d$. Now multiply both sides of the $y'_+$ version of
(\ref{eq3.51}) by $(1-\chi^\d)G(\dt,x)$ and integrate the result.
Integrate by parts and then use (\ref{eq3.25}) to find that
\begin{equation}
\label{eq3.56}
 y'_+(x)\leq \z_\d \exp(-\sqrt{r}/\z_\d) \int_X y'_+ .
\end{equation}
This last inequality with (\ref{eq3.29}) gives (\ref{eq3.49}) when
dist$(x,Z)\geq 2\d$. Thus, replacing $\d$ by $\d/2$ in (\ref{eq3.56}) gives
Proposition \ref{prop3.7}.

\section{The monotonicity formula}
\label{sec4}

Fix a ${\mbox{Spin}}^{\C}$ structure, a value of $r\geq 1$ and a solution
$(A,\psi)$ to the associated
version of (2.9). Let $B\subset X$ be an open set, and consider the
{\it energy} of $B$:
\setcounter{equation}{0}
\begin{equation}
\label{eq4.1}
 {\cal E}_B\equiv (4\sqrt{2})^{-1} r\int_B |\o| \
|(2^{-1/2}|\o| -|\psi|^2)| .
\end{equation}
Note that ${\cal E}_B\leq{\cal E}_X <\i$ by virtue of Lemma \ref{lem3.1} and
the second
line of (\ref{eq3.2}) in particular.  The purpose of this section is to first
estimate ${\cal E}_B$ from the above and from below in the case where $B$ is
a geodesic ball of some radius $s>0$. The second purpose will be to exploit
the estimates for the energy to refine some of the bounds in the previous
section.

\bs{\bf a) Monotonicity}

The following proposition describes the behavior of the energy ${\cal E}_B$
for the case where $B$ is
a geodesic ball in $X$ of some radius $s>0$.

\begin{proposition}
\label{prop4.1}
Fix a ${\mbox{Spin}}^{\C}$ structure for $X$. There is a constant $\z\geq
1$, and given $\d >0$, there
is a constant $\z_\d\geq 1$; and these constants have the following
significance: Fix $r\geq\z_\d$ and consider a solution $(A,\psi)$ to the
version of (2.9) which corresponds to the given ${\mbox{Spin}}^{\C}$
structure and $r$.
Let $B\subset X$ be a geodesic ball with center $x$ whose points all lie at
distance $\d$ or greater from $Z$. Let $s$ denote the radius of $B$ and
require $1/\z_\d \geq s\geq 2^{-1} r^{-1/2}$. Then:
\begin{itemize}
\item ${\cal E}_B \leq\z s^2$
\item ${\mbox{If}} \ \ |\a(x)| < (2\sqrt{2})^{-1}|\o|, \ \ {\mbox{then}} \ \
{\cal E}_B\geq\z^{-1}_\d s^2.$\eqnum
\end{itemize}
\end{proposition}

Proposition \ref{prop4.1} is proved in the next subsection. Note that this
proposition has the following crucial corollary:

\begin{lemma}
\label{lem4.2}
Fix a ${\mbox{Spin}}^{\C}$ structure of $X$. Given $\d >0$, there
is a constant $\z_\d >4$ with the following significance: Fix
$r\geq\z_\d$ and let $(A,\Psi)$ be a solution to (2.9) for the given
value of $r$ and the given ${\mbox{Spin}}^{\C}$ structure. Let $\rho\in
(\z_\d r^{-1/2},\z_\d^{-1}\d)$.
Then
\begin{itemize}
\item Let $\L$ be any set of disjoint balls of radius $\rho$ whose centers
lie on $\a^{-1}(0)$ and have distance at least $\d$ from $Z$. Then $\L$ has
less than $\z_\d\rho^{-2}$ elements.
\item The set of points in $\a^{-1}(0)$ with distance at least $\d$ from $Z$
has a cover by a set $\L$ of no more than $\z_\d\rho^{-2}$ balls of radius
$\rho$. Moreover, each ball in this set has center on $\a^{-1}(0)$ and
distance to $Z$ at least $\d/2$. Finally, the set of concentric balls of
radius $\rho/2$ is disjoint.
\end{itemize}
\end{lemma}

Note that Lemma \ref{lem4.2} plays the role in subsequent arguments that is
played
by Lemma $3.6$ in \cite{T1}.

\bs{\bf Proof of Lemma \ref{lem4.2}}\qua To prove the first assertion, use
Proposition \ref{prop4.1} to conclude that when $r$ is large, then the
energy of each
ball in the set $\L$ is at least $\z_\d^{-1}\rho^2$.
If there are $N$ such balls and they are all disjoint, then
${\cal E}_X\geq N\z_d^{-1}\rho^2$. Since ${\cal E}_X\leq\z$, this gives the
asserted bound on $N$. The second assertion follows from the first by
setting $\L'$ to equal a maximal (in number) set of disjoint balls of
radius $\rho/2$ whose centers lie on $\a^{-1}(0)$ and have distance at least
$\d/2$ from $Z$. With $\L'$ in hand, set $\L$ equal to the set whose balls
are concentric to those in $\L'$ but have radius $\rho$.

\bs{\bf b)  Proof of Proposition \ref{prop4.1}} The first two assertions
follow
from the following claim:  For fixed center $x$, consider ${\cal E}_B$ as a
function
of the radius $s$ of $B$.  Then ${\cal E}_B$ is a differentiable function
of $s$
which obeys the inequality:
\begin{equation}
\label{eq4.3}
{\cal E}_B \le 2^{-1}s(1 + \zeta_{\d}s)(1 + \zeta_{\d} r^{-1/2}) \frac
{d}{ds} {\cal
E}_B +
\zeta_{\d}s^4.
\end{equation}
If one is willing to accept (\ref{eq4.3}), then the proof of Proposition
\ref{prop4.1} proceeds by copying essentially verbatim that of Proposition
$3.1$ in
\cite{T1}.

With the preceding understood, the task at hand is to establish
(\ref{eq4.3}).  In
this regard, note that the argument for (\ref{eq4.3}) is only a slightly
modified
version of that for Proposition $3.2$ in \cite{T1}.  For this reason, the
discussion
below is brief.

To begin, remark that because of (2.9), one has
\begin{equation}
\label{eq4.4}
{\cal E}_B \le \int_B \o \wedge 2^{-1} i F_A.
\end{equation}
Meanwhile, $\o$ is exact on $B$, so can be written as $d\theta$ for some smooth
$1$--form on $B$.  Thus,
\begin{equation}
\label{eq4.5}
{\cal E}_B \le 2^{-1} \int_{\p B} \theta \wedge 2^{-1} i F_A.
\end{equation}

Since $\o$ is assumed to be nowhere vanishing on $B$, it follows that there
is a
coordinate system which is centered at $x$ and valid in a ball of radius
$\zeta_{\d}^{-1}$ about $x$ for which $\o$ pulls back to $\R^4$ as the
standard form
$\o_x = |\o(x)| \cdot (dy^1 \wedge dy^2 + dy^3 \wedge dy^4)$.  Moreover, this
coordinate chart can be chosen so that the pulled back metric is close to a
constant
multiple of the standard Euclidean metric on $\R^4$.  To be precise, one
can require
that the metric $g$ differ from $g_E = \sum_j dy^j \otimes dy^j$ as follows:
\begin{itemize}
\item $|g-g_E| \le \zeta_{\d}|y|$.
\item $|\p g| \le \zeta_{\d}$.\eqnum
\end{itemize}
%\begin{equation}
%\label{eq4.6}
%\end{equation}
Here, $\p g$ denotes the tensor of $y$--partial derivatives of $g$.  Note
that the
second line in (4.6) implies that the distance $s$ from the origin as
measured by the metric $g$ differs from that, $s_E$, measured by the
Euclidean metric
as follows:  $|s-s_E| \le \zeta_{\d}s^2$.

In these coordinates, the choice $\theta = 2^{-1}|\o(x)|(y^1dy^2 - y^2dy^1
+ y^3dy^4
- y^4dy^3)$ will be made.  Note that $|\theta|$ differs from $2^{-1}s$ by
no more
than $\zeta_{\d}s^2$.  With the preceding understood, it follows (as argued in
($3.21$--$24$) in \cite{T1}) that
\begin{equation}
\label{eq4.7}
{\cal E}_B \le 2^{-1}s(1 + \zeta_{\d}s)(1 + \zeta_{\d}r^{-1/2})4^{-1} r
\int_{\p B}
|\o(x)||(2^{-1}|\o| - |\a|^2)| + \zeta_{\d}s^4.
\end{equation}
Moreover, since $||\o| - |\o(x)|| \le \zeta_{\d}s|\o|$ on $\p B$, the
constant factor
$|\o(x)|$ above can be replaced by the variable factor $|\o|$ at the cost of
increasing $\zeta_{\d}$.  Thus, (\ref{eq4.7}) implies that
\begin{equation}
\label{eq4.8}
{\cal E}_B \le 2^{-1} s(1 + \zeta_{\d}s)(1 + \zeta_{\d}r^{-1/2})4^{-1}r
\int_{\p B}
|\o||(2^{-1}|\o| - |\a|^2)| + \zeta_{\d}s^4.
\end{equation}

To complete the argument for (\ref{eq4.3}), use Proposition \ref{prop3.3}
in the
previous section to replace the factor $(2^{-1}|\o| - |\a|^2)$ in
(\ref{eq4.8}) with
$(2^{-1}|\o| - |\psi|^2)$ at the cost of slightly increasing $\zeta_{\d}$.  The
resulting equation is (\ref{eq4.3}) after the identification of the $s$
derivative of
${\cal E}_B$ with the quantity $4^{-1} r\int_{\p B} |\o|(2^{-1}|\o| -
|\psi|^2)$.

\bs{\bf a) A refined curvature bound}

The results in Proposition \ref{prop4.1} about ${\cal E}_B$ can be used to
refine the
bound in Proposition \ref{prop3.4} for $|F_A^-|$.  The following proposition
summarizes:

\begin{proposition}
\label{prop4.3}
Fix a $\mbox{Spin}^{\C}$ structure $s$, and fix $\d > 0$.  Then, there exist
constants $\zeta_{\d}$, $\zeta'_{\d} \ge 1$ with the following
significance:  Let
$(A,\psi)$ be a solution to the $s$ version of (2.9) where $r \ge
\zeta_{\d}$.  Then, at points in $X$ with distance $\d$ or more from $Z$,
\begin{equation}
\label{eq4.9}
|F_A^-| \le r(2\sqrt{2})^{-1}(2^{-1/2}|\o| - |\a|^2) + \zeta'_{\d}.
\end{equation}
\end{proposition}

The remainder of this section is occupied with the

\bs{\bf Proof of Proposition \ref{prop4.3}}\qua The proof amounts to a slight
modification of the arguments which prove Proposition $3.4$ of \cite{T1}.
To start,
introduce the functions $q_2$ and $q_3$ as in (\ref{eq3.38}) and
(\ref{eq3.39}).
Because of (\ref{eq3.40}) and (\ref{eq3.41}) one has $|q_2| \le
\zeta_{\d}$, and so
(due to Proposition \ref{prop3.3}'s bound on $|\b|^2$), it is enough to
bound $q_3$
by a uniform constant.  In this regard, note that $q_3$ obeys the equation
\begin{equation}
\label{eq4.10}
2^{-1}d^*dq_3 + 4^{-1}|\a|^2 q_3 \le \zeta_{\d}r.
\end{equation}
where $\s \ge \d$ and when $r \ge \zeta_{\d}$.  Also, $q_3 = 0$ where $\s =
\d$.

With these last points understood, the key to the proof is the following lemma
(compare with Lemma $3.5$ in \cite{T1}):

\begin{lemma}
\label{lem4.4}
Fix a $\mbox{Spin}^{\C}$ structure $s$, and fix $\d > 0$.  Then, there is a
constant
$\zeta_{\d} \ge 1$ with the following significance:  Let $(A,\psi)$ be a
solution to
the $s$ version of (2.9) where $r \ge \zeta_{\d}$.  Then, there is
a smooth
function $u$ which is defined on the set of points in $X$ with distance
$\d$ or more
from $Z$ and which obeys
\begin{itemize}
\item $|u| \le \zeta_{\d}$.
\item $2^{-1} d^*d u \ge r$ where $|\a| \le (2\sqrt{2})^{-1} |\o|$.
\item $|d^*d u| \le \zeta_{\d}r$.
\item $u = 0$ where $\s = \d$.{\rm\eqnum}
\end{itemize}
%\begin{equation}
%\label{eq4.11}
%\end{equation}
\end{lemma}

The proof of Proposition \ref{prop4.3} given Lemma \ref{lem4.4} is
essentially the
same as that of Proposition $3.4$ in \cite{T1} given Lemma $3.5$ in \cite{T1}.
Meanwhile, the proof of Lemma \ref{lem4.2} is a Dirichlet boundary
condition version
of the proof of Lemma $3.5$ in \cite{T1}.  The modifications to the
argument for the
latter in \cite{T1} are straightforward and left to the reader.

\section{Local properties of $\a^{-1}(0)$}
\label{sec5}

The purpose of this section is to summarize some of the local properties of
$\a^{-1}(0)$ at points in the complement of $Z$.  At issue here is the
behavior of
$\a^{-1}(0)$ at length scales of order $r^{-1/2}$.

The strategy for the investigation at such scales is as follows:  Fix $\d >
0$ and a
point $x$ whose distance from $Z$ is at least $\d$.  A Gaussian coordinate
system
based at $x$ defines an embedding $h\co  \R^4 \rightarrow X$ which maps the
origin to
$x$ and which sends straight lines through the origin in $\R^4$ to
geodesics in $X$
through $x$.  Moreover, the pull-back via $h$ of the Riemannian metric
agrees with
the Euclidean metric to second order at the origin.  The Gaussian
coordinate charts
at $x$ are parametrized by the group $SO(4)$ (to be precise, the fiber of
the frame
bundle at $X$).  In particular, there are Gaussian coordinate systems at
$X$ which
pull $\o$ back as $h^*\o = |\o(x)|(dy^1 \wedge dy^2 + dy^3 \wedge dy^4) + {\cal
O}(|y|)$.  Such a Gaussian coordinate system will be called a {\em complex
Gaussian
coordinate system}.  Indeed, a Gaussian coordinate system at $x$ is called
complex
precisely when the differential of the corresponding $h$ at the origin
intertwines
the standard almost complex structure on $\R^4 = \C^2$ with $J|_x$.  The
complex
Gaussian coordinate systems at $x$ are parametrized by the $U(2)$ subgroup
of $SO(4)$.

Now, fix a $\mbox{Spin}^{\C}$ structure on $X$ and $r \ge 1$ and then let
$(A,\psi)$
be a solution to the corresponding version of (2.9).  Then,
pull-back by the
map $h$ of a Gaussian coordinate system at some $x$ defines $(A,\psi)$ as
fields open
$\R^4$.

Given $\l > 0$, define the dilation map $\d_{\l}\co  \R^4 \rightarrow \R^4$ by its
action on the coordinate functions $y\co  \d_{\l}^* y = \l^{-1}y$.  With $x$
chosen in
the complement of $Z$ set $\l = (r|\o(x)|)^{1/2}$ and let $h$ be a complex
Gaussian
coordinate system based at $x$.  Given $(A,\psi = (\a,\b))$, define the data
$(\underline{A},(\underline{\a},\underline{\b}))$ on $\R^4$ by the rule:
\setcounter{equation}{0}
\begin{equation}
\label{eq5.1}
(\underline{A},(\underline{\a},\underline{\b})) \equiv
\d_{\l}^*h^*(A,|\o(x)|^{-1/2}(\a,\b)).
\end{equation}

The plan now is to compare
$(\underline{A},(\underline{\a},\underline{\b}))$ with
some standard objects on $\R^4$.  These standard objects are discussed in
Proposition
$4.1$ of \cite{T1}.  The following digression constitutes a brief summary:  A
connection $a_0$ on the trivial complex line bundle over $\R^4$ and a
section $\a_0$
of this line bundle will be said to a {\em solution to the Seiberg--Witten
equations}
on $\R^4$ when the following conditions hold:

\begin{itemize}
\item The curvature $2$--form, $F_a$, of $a_0$ is of type $1\!\!-\!\!1$ 
with respect to the standard almost complex structure on $\R^4$ and so 
defines a holomorphic structure
(and associated ${\bar \p}$ operator) on the trivial complex line bundle.
\item The section $\a_0$ is holomorphic with respect to the $a_0$--complex
structure on
the trivial complex line bundle.
\item $F_a^+ = -i8^{-1}(1 - |\a_0|^2)(dy^1 \wedge dy^2 + dy^3 \wedge dy^4)$
\item $|\a_0| \le 1$.
\item $|F_a^-| \le |F_a^+| \le (4\sqrt{2})^{-1}(1-|\a_0|^2)$.
\item $|\nabla_a\a_0| \le z(1-|\a_0|^2)$.
\item For each $N \ge 1$, the integral of $(1-|\a_0|^2)$ over the ball of
radius $N$
is bounded by $zN^2$.\eqnum
\end{itemize}
%\begin{equation}
%\label{eq5.2}
%\end{equation}
Here, $z$ is a constant.  (Note that these conditions differ from the
conditions
listed in $(4.3)$ of \cite{T1} in that no assumption on the integrability of
$|F_a^+|^2 - |F_a^-|^2$ is made here.  It is most probably true that the latter
condition is a consequence of those in (5.2).)

The following proposition summarizes the basic properties of solutions to
(5.2):

\begin{proposition}
\label{prop5.1}
Let $(a_0,\a_0)$ obey the conditions in (5.2).  Then:
\begin{itemize}
\item Either $|\a_0| < 1$ everywhere or else $|\a_0| = 1$ and
$(a_0,\a_0)$ is gauge equivalent to the trivial solution $(a_0 =
0,\a_0 = 1)$.  In the former case, 
$\a_0^{-1}(0) \ne \emptyset$ and $\a_0^{-1}(0)$ is the zero set of a
polynomial in
the complex coordinates for $\R^4 = \C^2$.
\item Either $|F_a^-| < |F_a^+|$ everywhere or else $|F_a^-| \equiv
|F_a^+|$ and
there is a $\C$--linear map $s\co  \C^2 \rightarrow \C$ and a solution
$(a_1,\a_1)$ to
the {\em vortex equations} on $\C$ with the property that $(a_0,\a_0)$ is gauge
equivalent to the pull-back $s^*(a_1,\a_1)$.  In this case, $\a_0^{-1}(0)$ is a
finite set of parallel, complex planes.
\item Given the constant $z$ in (5.2), there is an upper bound on
the order
of vanishing of $\a_0$ at any point in $\C^2$.
\item The set of gauge equivalence classes of $(a_0,\a_0)$ which obey
(5.2)
for a fixed value of $z$ is sequentially compact with respect to convergence on
compact subsets of $\R^4$ in the $C^{\i}$ topology.
\item Given the value of $z$ in (5.2), there exists $z_1 > 0$ such that
\begin{equation}
\label{eq5.3}
(1-|\a_0|^2) + |\nabla_a\a_0|^2 \le z_1
\exp[-\mbox{\rm dist}(\cdot,\a_0^{-1}(0))/z_1].
\end{equation}
\end{itemize}
\end{proposition}

This proposition restates various assertions of Proposition $4.1$ in
\cite{T1}; the
reader is referred to Section~4e of \cite{T1} for the proof.

End the digression.

The relevance of the solutions to the standard Seiberg--Witten solutions to the
problem at hand is summarized by the next proposition:

\begin{proposition}
\label{prop5.2}
Fix a $\mbox{Spin}^{\C}$ structure for $X$.  Given $\d > 0$, there is a
constant
$z_{\d} \ge 1$, and given $R \ge 1$, $k \ge 1$ and $\e > 0$, there is another
constant $\zeta_{\d}$ and these have the following significance:  Let $r >
\zeta_{\d}$ and let $(A,\psi)$ be a solution to (2.9) as defined
with the
given $\mbox{Spin}^{\C}$ structure and with $r$.  Suppose that $x \in X$
has distance
at least $\d$ from $Z$.  Now define the fields
$(\underline{A},(\underline{\a},\underline{\b}))$ as in (\ref{eq5.1}).
Then there
exists a solution $(a_0,\a_0)$ to the $z_1 \le z_{\d}$ version of
(5.2) and a
gauge transformation $\varphi\co  \C^2 \rightarrow S^1$ such that
$\varphi^*(\underline{A},(\underline{\a},\underline{\b})) -
(2a_0,(\a_0,0))$ has
$C^k$ norm less than $\e$ in the ball of radius $R$ and center $0$ in $\R^4$.
Furthermore, $|\a_0|$ is not constant if $|\a(x)| \le 1 - \e$.
\end{proposition}

{\bf Proof of Proposition \ref{prop5.2}}\qua The proof of this
proposition can be
found by lifting from Section~4c of \cite{T1} the proofs of the analogous
assertions
of Proposition $4.2$ of \cite{T1}.  (The lack of control here on the
integral over
$\R^4$ of $|F_a^+|^2 - |F_a^-|^2$ precludes only the use of the proofs in
\cite{T1}
of statements which actually refer to this integral.)

\section{Large $r$ behavior away from $\a^{-1}(0)$}
\label{sec6}

Fix a $\mbox{Spin}^{\C}$ structure for $X$ and then consider a solution
$(A,\psi)$ to
(2.9) for the given $\mbox{Spin}^{\C}$ structure and for some large
value of
$r$.  The purpose of this section is to investigate the behavior of
$(A,\psi)$ to
(2.9) at points which lie neither on $Z$ nor on $\a^{-1}(0)$.  Here
are the
basic observations:  First, $2^{-1/2}|\o| - |\a|^2$ and $|\b|$ are both ${\cal
O}(r^{-1})$.  In particular, this means that $F_A$ is bounded.  More to the
point,
the connection $A$ is close to a canonical connection $A^0$ whose gauge
orbit depends
only on the metric and the choice of $\o$.  (Note that this orbit is
independent of
the chosen $\mbox{Spin}^{\C}$ structure.)  Proposition \ref{prop6.1},
below, gives
the precise measure of closeness that is used here.

The statement of Proposition \ref{prop6.1} requires a preliminary, three part
digression whose purpose is to define the connection $A^0$.  The first part
of the
digression remarks that the $-\sqrt{2}i|\o|$ eigenspace of the Clifford multiplication
endomorphism by $c_+(\o)$ on $S_+$ defines a complex line bundle $E
\rightarrow X -
Z$.  The component $\a$ of $\psi$ is a section of $E$, and then the
component $\b$ is
one of $K^{-1}E$.  Here, $K^{-1}$ is the inverse of the canonical bundle,
$K$, for
the almost complex structure $J \equiv \sqrt{2} g^{-1}\o/|\o|$ on $X - Z$.

Note that the line bundle $L = \mbox{det}(S_+)$ restricts to $X-Z$ as
\setcounter{equation}{0}
\begin{equation}
\label{eq6.1}
L|_{X-Z} \approx K^{-1}E^2.
\end{equation}
Meanwhile, $\a$ trivializes $E$ where $\a \ne 0$ and so the unit length section
$\a^2/|\a|^2$ of $E^2$ provides an isometric identification of $L$ with
$K^{-1}$ on
the complement of $Z$ and $\a^{-1}(0)$.

The second part of the digression reviews the definition from \cite{T4} or
Section~1c
of \cite{T1} of a canonical connection on the line bundle $K^{-1} \rightarrow
X - Z$.
(Note that this connection is unique up to gauge equivalence.)  To define a
canonical
connection, first remark that there is a canonical $\mbox{Spin}^{\C}$
structure for
$X-Z$ so that the $-\sqrt{2}i|\o|$ eigenbundle for the $c_+(\o)$ action on the
corresponding
$S_+$ is the trivial bundle over $X-Z$.  For this $\mbox{Spin}^{\C}$
structure, the
corresponding line bundle $L$ is isomorphic to $K^{-1}$.  With this
understood, there
is a unique connection (up to isomorphism) on $K^{-1}$ for which the induced
connection on the aforementioned $-\sqrt{2}i|\o|$ eigenbundle is trivial.  Such a
connection
is a canonical one.

Part~3 of the digression defines a canonical connection on $L|_{X-Z}$ by
using the
identification $\a^2/|\a|^2$ between $L$ and $K^{-1}$ (where $\a \ne 0$ on
$X-Z$) to
pull a canonical connection on $K^{-1}$ back to $L$.

End the digression.

\begin{proposition}
\label{prop6.1}
Fix a $\mbox{Spin}^{\C}$ structure for $X$ and $\d > 0$.  There is a constant
$\zeta_{\d} \ge 1$ with the following significance:  Suppose that $r \ge
\zeta_{\d}$
and that $(A,\psi)$ are a solution to (2.9) as defined by $r$ and
the given
$\mbox{Spin}^{\C}$ structure.  There is a canonical connection $A^0$ on
$L|_{X-Z}$
for which $|A-A^0| + |F_A-F_{A^0}| \le \zeta_{\d}r^{-1} + \zeta_{\d} r
\exp[-\sqrt{r}
\mbox{\rm\ dist}(x,\a^{-1}(0))/\zeta_{\d}]$ at all points $x \in X$ with
distance $\d$ or
more from $Z$ and distance $r^{-1/2}$ or more from $\a^{-1}(0)$.
\end{proposition}

(This proposition should be compared with Proposition $4.4$ in
\cite{T1}.)

\bs{\bf Proof of Proposition \ref{prop6.1}}\qua Away from $Z$ and where $\a
\ne 0$,
the section $\a/|\a|$ defines a trivialization of the line bundle $E$, and
with this
understood, the difference between $A$ and a particular canonical
connection $A^0$ is
given by $2({\bar \a}/|\a|)\nabla_A(\a/|\a|)$.  Thus, the absolute value of
$\nabla_A(\a/|\a|)$ measures the size of $A - A^0$.  Likewise, the norm of $d_A
\nabla_A(\a/|\a|)$ measures the size of $F_A - F_{A^0}$.

With the task ahead now clear, note that the arguments which establish the
required
bounds on the derivatives of $\a$ are, for the most part, straightforward
modifications of the arguments which prove Proposition $4.4$ in \cite{T1}.  In
particular, the reader will be referred to the latter reference at numerous
points.
In any event, the details are given in the subsequent four steps.

\medskip
{\bf Step 1}\qua A straightfoward modification of the proof of
Proposition $4.4$ in \cite{T1} (which is left to the reader) proves the
following
preliminary estimate:

\begin{lemma}
\label{lem6.2}
Fix a $\mbox{Spin}^{\C}$ structure for $X$ and fix $\d > 0$.  There is a
constant
$\zeta_{\d}$ with the following significance:  Let $(A,\psi)$ solve the
version of
(2.9) which is defined by the given $\mbox{Spin}^{\C}$ structure
and by $r
\ge \zeta_{\d}$.  If $x \in X$ has distance $\d$ or more from $Z$, then
\begin{equation}
\label{eq6.3}\begin{split}
r|(2^{-1/2}|\o|-|\a|^2)|7&+ r^2|\b|^2 + |\nabla_A\a|^2 + r|\nabla_A\b|^2\\
& \le
\zeta_{\d}(1 + r \exp[-\sqrt{r} \mbox{ dist}(x,\a^{-1}(0))/\zeta_{\d}]).
\end{split}\end{equation}
\end{lemma}

{\bf Step 2}\qua Now, introduce $\underline{\a} \equiv
2^{1/4}|\o|^{-1/2}\a$.  Add the two lines of (\ref{eq3.53}) to obtain a differential
inequality for the function $y \equiv |\nabla_A\underline{\a}|^2 +
|\nabla_A\beta|^2$.  Use Lemma \ref{lem6.2} to bound $|w|$ to find that the
aforementioned inequality implies that $2^{-1}d^*d y +
4^{-1}r\zeta_{\d}^{-1} y \le
\zeta_{\d}r^{-1}$ at points with distance $\d$ or more to $Z$ and distance
$\zeta_{\d}r^{-1/2}$ or more to $\a^{-1}(0)$.  Note that (\ref{eq6.4})
implies that
$y' \equiv y - 4\zeta_{\d}^2 r^{-2}$ obeys the inequality
\begin{equation}
\label{eq6.4}
2^{-1}d^*dy' + 4^{-1}r\zeta_{\d}^{-1}y' \le 0.
\end{equation}
at points with distance $\d$ or more to $Z$ and $\zeta_{\d}r^{-1/2}$ or more to
$\a^{-1}(0)$.

Given (\ref{eq6.4}), a straightforward modification of the proof of
Proposition $4.4$
in \cite{T1} yields the bound
\begin{equation}
\label{eq6.5}
|\nabla_A\underline{\a}|^2 + |\nabla_A\beta|^2 \le \zeta_{\d}r^{-2} +
\zeta_{\d} r
\exp(-\sqrt{r} \mbox{ dist}(\cdot,\a^{-1}(0))/\zeta_{\d})
\end{equation}
at points with distance $2\d$ or more from $Z$.  (Bounds on the size of both
$|\nabla_A\underline{\a}|^2$ and $|\nabla_A\beta|^2$ near $\a^{-1}(0)$ come via
Proposition \ref{prop3.7}.)

Take the $\d/2$ version of (\ref{eq6.5}) with the fact that
$\underline{\a}/|\underline{\a}| = \a/|\a|$ to bound the difference between
$A$ and a
canonical connection on $L$ by $$\zeta_{\d}(r^{-1} + r^{1/2} \exp(-\sqrt{r}
\mbox{ dist}(\cdot,\a^{-1}(0))/\zeta_{\d}))$$ at points with distance $\d$
or more
from $Z$ and $r^{-1/2}$ or more from $\a^{-1}(0)$.

\medskip
{\bf Step 3}\qua As remarked above, a bound on
$|\nabla_A^2\underline{\a}|$ provides a bound on $|F_A - F_{A^0}|$.  To
obtain the
latter, first differentiate (\ref{eq3.55}) and commute covariant derivatives to
obtain an equation for $\nabla_A^2\underline{\a}$ of the form
$\nabla_A^*\nabla_A(\nabla_A^2\underline{\a}) +
(4\sqrt{2})^{-1}r|\o|(\nabla_A^2\underline{\a}) + \mbox{ Remainder} = 0$.
Take the
inner product of this last equation with $\nabla_A^2\underline{\a}$ to
obtain an
equation for $|\nabla_A^2\underline{\a}|^2$ having the form $2^{-1}
d^*d|\nabla_A^2\underline{\a}|^2 + 4^{-1} r|\o||\nabla_A^2\underline{\a}|^2 +
|\nabla_A(\nabla_A^2\underline{\a})|^2 +
\<\nabla_A^2\underline{\a},\mbox{Remainder}\> = 0$.  Here, $\<\ ,\ \>$
denotes the
Hermitian inner product on $E \otimes (\otimes_2 T^*X)$.  A similar
equation for
$|\nabla_A^2\b|$ is obtained by differentiating (\ref{eq3.56}).  Add the
resulting
two equations.  Then, judicious use Lemma \ref{lem6.2}, (\ref{eq6.6}) and the
triangle inequality produces a differential inequality for $y' \equiv
|\nabla_A^2\underline{\a}|^2 + |\nabla_A^2\b|^2 - \zeta_{\d}r^{-2}$ which
has the
same form as (\ref{eq6.4}).  And, with this understood, the arguments which
yield
(\ref{eq6.5}) yield the bound
\begin{equation}
\label{eq6.6}
y' \le \zeta_{\d}r^{-2} + \zeta_{\d}(\sup_{\s \ge \d} |y'|)\exp(-\sqrt{r}
\mbox{
dist}(\cdot,\a^{-1}(0))/\zeta_{\d})
\end{equation}
at points where $\s \ge 2\d$.

\medskip
{\bf Step 4}\qua The $\d/2$ version of (\ref{eq6.6}) with a
bound on
$|\nabla_A^2\underline{\a}|^2 + |\nabla_A^2\b|^2$ where $\s \ge \d/2$ by
$\zeta_{\d}r^2$ gives Proposition \ref{prop6.1}'s bound on $|F_A-F_{A^0}|$.
Thus,
the last task is to obtain a supremum bound on $|\nabla_A^2\underline{\a}|^2 +
|\nabla_A^2\b|^2$.

For this purpose, fix a ball of radius $2r^{-1/2}$ whose points all have
distance
$\d/4$ or more from $Z$.  Take Gaussian coordinates based at the center of
this ball
and rescale so that the radius $r^{-1/2}$ concentric ball becomes the
radius $1$ ball
in $\R^4$ with center at the origin.  Equation (2.9) rescales to
give an
$r=1$ version of the same equation on the radius $2$ ball in $\R^4$ with a
metric
$g'$ which is close to the Euclidean metric $g_E$ and form $\o'$ which is
close a
constant self dual form of size $|\o(x)|$.  Here, $x$ is the center of the
chosen
ball in $X$.  To be precise, $|g' - g_E| \le \zeta r^{-1}$ and the
derivatives of
$g'$ of order $k \ge 2$ are ${\cal O}(r^{-k/2})$ in size.  Meanwhile the
form $\o'$
differs by ${\cal O}(r^{-1/2})$ from a constant form, and its $k$-th
derivatives are
${\cal O}(r^{-k/2})$ in size.

With the preceding understood, standard elliptic regularity results (as in
Chapter~6
of \cite{My}) bound the second derivatives of the rescaled versions of
$\underline{\a}$ and $\b$ by $\zeta_{\d}$.  Rescaling the latter bounds
back to the
original size gives $|\nabla_A^2\underline{\a}|^2 + |\nabla_A^2\b|^2 \le
\zeta_{\d}r^2$ as required.

\section{Proof of Theorem \ref{th2.2}}
\label{sec7}

Fix a $\mbox{Spin}^{\C}$ structure for $X$ and suppose that there exists an
unbounded, increasing sequence $\{r_n\}$ of positive numbers with the
property that
each $r = r_n$ version of (2.9) with the given $\mbox{Spin}^{\C}$
structure
has a solution $(A_n,\psi_n)$.  The purpose of this section is to
investigate the $n
\rightarrow \i$ limits of $\a_n^{-1}(0)$ and in doing so, prove the claims
of Theorem
\ref{th2.2}.  This investigation is broken into six parts.

\bs{\bf a)  The curvature as a current}

Each connection $A_n$ has its associated curvature $2$--form, and the difference
between $A_n$'s curvature $2$--form and the curvature $2$--form of the canonical
connection on $K^{-1}$ will be viewed as a current on $X$.  This current,
${\cal
F}_n$, associates to a smooth $2$--form $\mu$ the number
\setcounter{equation}{0}
\begin{equation}
\label{eq7.1}
{\cal F}_n(\mu) \equiv 2^{-1} \int_X \frac {i}{2\pi} (F_{A_n} - F_{A^0})
\wedge \mu.
\end{equation}
Here, $F_{A^0}$ is the curvature $2$--form of the canonical connection on
$K^{-1}$.
(Even though the canonical connection on $K^{-1}$ is defined only over
$X-Z$ (see the
beginning of the previous section), the norm of its curvature is none-the-less
integrable over $X$.  Thus, (\ref{eq7.1}) makes sense even for $\mu$ whose
support
intersects $Z$.  The integrability of $|F_{A^0}|$ follows from the bound
$|F_{A^0}|
\le \zeta \mbox{ dist}(\cdot,Z)^{-2}$.)

With the sequence $\{{\cal F}_n\}$ understood, fix $\d > 0$ and suppose
that each
point in the support of $\mu$ has distance $\d$ or more from $Z$.  It then
follows
from Lemma \ref{lem3.1} and Proposition \ref{prop4.3} that
\begin{equation}
\label{eq7.2}
|{\cal F}_n(\mu)| \le \zeta_{\d} \sup_X |\mu|.
\end{equation}
This uniform bound implies that the sequence $\{{\cal F}_n\}$ defines a bounded
sequence of linear functional on the space of smooth $2$--forms on $X$ with
support
where the distance to $Z$ is at least $\d$.

With the preceding understood, a standard weak convergence argument finds a
subsequence of $\{{\cal F}_n\}$ (hence renumbered consecutively) which
converges in
the following sense:  Let $\mu$ be a smooth $2$--form with compact support
on $X-Z$
and then $\lim_{n \rightarrow \i} {\cal F}_n(\mu)$ exists.  Moreover, this
limit,
\begin{equation}
\label{eq7.3}
{\cal F}(\cdot) \equiv \lim_{n \rightarrow \i} {\cal F}_n(\cdot),
\end{equation}
defines a bounded linear functional when restricted to forms whose support has
distance from $Z$ which is bounded from below by any fixed positive number.

Note that the current ${\cal F}$ is {\em integral} in the following sense:
Let $\mu$
be a closed $2$--form with compact support on $X-Z$ and with integral periods on
$H_2(X-Z;\Z)$.  Then
\begin{equation}
\label{eq7.4}
{\cal F}(\mu) \in \Z.
\end{equation}

\bs{\bf b)  The support of ${\cal F}$}

This part of the discussion considers the support of the current ${\cal
F}$.  Here is
the crucial lemma:

\begin{lemma}
\label{lem7.1}
There is a closed subspace $C \subset X-Z$ with the following properties:

\begin{itemize}
\item ${\cal F}(\mu) = 0$ when $\mu$ is a $2$--form on $X$ with compact
support in
$(X-Z) - C$.
\item Conversely, let $B \subset X-Z$ be an open set which intersects $C$.
Then
there is a $2$--form $\mu$ with compact support in $B$ and with ${\cal
F}(\mu) \ne 0$.
\item Fix $\d > 0$.  Then the set of points in $C$ with distance at least
$\d$ from
$Z$ has finite $2$--dimensional Hausdorff measure.
\item Conversely, let $\d > 0$ and there is a constant $\zeta_{\d} \ge 1$
with the
following significance:  Let $\rho \in (0,\zeta_{\d}^{-1})$ and let $B
\subset X$ be
a ball of radius $\rho$ and center on $Z$ whose points have distance $\d$
or more
from $Z$.  Then the $2$--dimensional Hausdorff measure of $B \cap C$ is
greater than
$\zeta_{\d}^{-1}\rho^2$.
\item There is a subsequence of $(A_n,\psi_n)$ such that the corresponding
sequence
$\{\a_n^{-1}(0)\}$ converges to $C$ in the following sense:  For any $\d >
0$, the
following limit exists and is zero:
\begin{equation}
\label{eq7.5}
\lim_{n \rightarrow \i} [\sup_{\{x\in C: \mbox{\footnotesize\rm\ dist}(x,Z) 
\ge \d\}} \mbox{\rm
dist}(x,\a_n^{-1}(0)) + \sup_{\{x \in \a_n^{-1}(0): \mbox{\footnotesize\rm\ 
dist}(x,Z) \ge \d\}}
\mbox{\rm dist}(x,C))].
\end{equation}
\end{itemize}
\end{lemma}

{\bf Proof of Lemma \ref{lem7.1}}\qua To construct $C$, consider first a
large
positive integer $N$ and a very large positive integer $n$.  (Here, a lower
bound on
$n$ comes from the choice of $N$.)  Use Lemma \ref{lem4.2} to find a set
$\L'_n(N)$
of balls of radius $16^{-N}$ with the following properties:  The balls are
disjoint,
their centers lie on $\a_n^{-1}(0)$ and have distance at least $4 \cdot
16^{-N}$ from
$Z$, and the set $\L_n(N)$ of concentric balls of radius $2 \cdot 16^{-N}$
covers the
set of point in $\a_n^{-1}(0)$ with distance $8 \cdot 16^{-N}$ from $Z$.
According
to Lemma \ref{lem4.2}, when $n$ is sufficiently large, this set $\L_n(N)$
has a bound
on the number of its elements which is independent of $n$.  Let $\nu(N)$ be an
integer which is greater than the number of elements in each $\L_n(N)$.

Label the centers of the balls in $\L_n(N)$ and then add the final point
some number
of times (if necessary) to make a point $\underline{x}_n(N) \equiv
(x_n(1;N),\dots,x_n(\nu(N);N))$ $\in \times_{\nu(N)}X$.  By a diagonalization
process, one
can find an infinite sequence of indices $n$ (hence relabeled consecutively
from $1$)
so that for each $N$, the sequence $\{\underline{x}_n(N)\}_{n \ge 1}$
converges in
$\mathbf{\times}_{\nu(N)} X$.

For each $N$, let $\underline{x}(N) \equiv \{x(1,N),\dots,x(\nu(N),N)\}$
denote the
limit of\break $\{\underline{x}_n(N)\}_{n \ge 1}$.  One can think of
$\underline{x}(N)$ as
either a point in $\times_{\nu(N)} X$, or else an ordered set of $\nu(N)$ points
in $X$.
Think of $\underline{x}(N)$ in the latter sense, and let $U(N)$ denote the
union of
the balls of radius $4 \cdot 16^{-N}$ with centers at the points
$\underline{x}(N)$
(that is, at $\{x(i,N)\}_{1 \le i \le \nu(N)}$).  Lemma $5.1$ of \cite{T1}
argues
that these sets are nested in that $U(N+1) \subset U(N)$.  With this
understood, set
\begin{equation}
\label{eq7.6}
C \equiv \bigcap_{N \ge 1} U(N).
\end{equation}

The argument for the asserted properties of $C$ is essentially the same as
that for
Lemma $5.2$ in Section~5c of \cite{T1}.  In this regard, Lemma \ref{lem4.2}
here
replaces Lemma $3.6$ in \cite{T1}, and Proposition \ref{prop4.3} here replaces
Proposition $4.4$ in \cite{T1}.  The details are straightforward and left
to the
reader.

\bs{\bf c)  A positive cohomology assignment}

The purpose of this subsection is to give a more precise characterization
of the
distribution ${\cal F}$.

To begin, note that the construction of $C$ indicated above can be used (as
in the
proof of Lemma $5.3$ in \cite{T1}) to prove that the current ${\cal F}$ is
{\em type}
$1\!\!-\!\!1$ in the sense that ${\cal F}(\mu) = 0$ when $\mu$ is a section of the
subbundle
$K^{-1} \subset \L_+$.  The fact that ${\cal F}$ is type $1\!\!-\!\!1$ is implied
by Lemma
\ref{lem7.2} below which is a significantly stronger assertion:

\begin{lemma}
\label{lem7.2}
Let $D \subset \C$ denote the standard, unit disc, and let $\s \co  D
\rightarrow X
- Z$ be a smooth map which extends to the closure, $\underline{D}$, of $D$ as a
continuous map that sends $\p \underline{D}$ into $X-C$.  Then,
$\left\{2^{-1} \int_D
\s^*\left( \frac {i}{2\pi} (F_{A_n} - F_{A^0})\right)\right\}_{n \ge 1}$
converges,
and the limit, $I(\s) \in \Z$. Moreover,

\begin{itemize}
\item $I(\s) = 0$ if $\s(D) \cap C = \emptyset$.
\item $I(\s) > 0$ if $\s$ is a pseudo-holomorphic map and $\s^{-1}(C) \ne
\emptyset$.{\rm\eqnum}
\end{itemize}
%\begin{equation}
%\label{eq7.7}
%\end{equation}
\end{lemma}

{\bf Proof of Lemma \ref{lem7.2}}\qua The arguments which prove Lemma
$6.2$ in
\cite{T1} are of a local character and so can be brought to bear directly
to give
Lemma \ref{lem7.2}.

The discussion surrounding Lemma $6.2$ in \cite{T1} concerns the notion of
a {\em
positive cohomology assignment} for $C$.  The latter is defined as follows:
First,
let $D \subset \C$ be the standard unit disk again.  A map $\s \co  D
\rightarrow X
- Z$ is called {\em admissible} when $\s$ extends as a continuous map to
the closure,
$\underline{D}$, of $D$ which maps $\p \underline{D}$ into $X-C$.  A positive
cohomology assignment specifies an integer, $I(\s)$, for each admissible
map $\s$ from
$D$ to $X-Z$ subject to the following constraints:

\begin{itemize}
\item If $\s(D) \subset X-C$, then $I(\s) = 0$.
\item A homotopy $h\co  [0,1] \x D \rightarrow X$ is called admissible when it
extends as
a continuous map from $[0,1] \x \underline{D}$ into $X$ that sends $[0,1] \x
\p\underline{D}$ to $X-C$.  If $h$ is an admissible homotopy, then
$I(h(1,\cdot)) =
I(h(0,\cdot))$.
\item Let $\s \co  D \rightarrow X$ be admissible and suppose that $\theta\co  D
\rightarrow D$ is a proper, degree $k$ map.  Then $I(\s\cdot\theta) = kI(\s)$.
\item Suppose that $\s \co  D \rightarrow X$ is admissible and that
$\s^{-1}(C)$ is
contained in a disjoint union $\cup_{\nu}D_{\nu} \subset D$, where each
$D_{\nu} =
\theta_{\nu}(D)$ with $\theta_{\nu}\co  D \rightarrow D$ being an orientation
preserving
embedding.  Then $I(\s) = \sum_{\nu} I(\s\cdot \theta_{\nu})$.
\item If $\s$ is admissible and pseudo-holomorphic with $\s^{-1}(C) \ne
\emptyset$,
then $I(\s) > 0$.\eqnum
\end{itemize}
%\begin{equation}
%\label{eq7.8}
%\end{equation}

The next result follows from Lemma \ref{lem7.2} and the particulars of the
definition
of ${\cal F}$ in (\ref{eq7.3}) as a limit.

\begin{lemma}
\label{lem7.3}
Let $C$ be as in Lemma \ref{lem7.1} and let $I(\cdot)$ be as described in Lemma
\ref{lem7.2}.  Then $I(\cdot)$ defines a positive cohomology assignment for
$C$.
\end{lemma}

{\bf Proof of Lemma \ref{lem7.3}}\qua Use the proof of Lemma $6.2$ in
\cite{T1}.

\bs{\bf d) $C$ as a pseudo-holomorphic submanifold}

Proposition $6.1$ in \cite{T1} asserts that a closed set in a compact,
symplectic
$4$--manifold with finite $2$--dimensional Hausdorff measure and a positive
cohomology
assignment is the image of a compact, complex curve by a pseudo-holomorphic
map.
(Note, however that the proof of Proposition $6.1$ in \cite{T1} has errors
which
occur in Section~6e of \cite{T1}, and so the reader is referred to the
revised proof
in the version which is reprinted in \cite{T6}.)  It is important to
realize that the
assumed compactness of $X$ in the statement of Proposition $6.1$ of
\cite{T1} is
present only to insure that the complex curve in question is compact.  In
particular,
the proof of Proposition $6.1$ in \cite{T1} from the reprinted version in
\cite{T6}
yields:

\begin{proposition}
\label{prop7.4}
Let $Y$ be a $4$--dimensional symplectic manifold with compatible almost complex
structure.  Suppose that $C \subset Y$ is a closed subset with the following
properties:

\begin{itemize}
\item The restriction of $C$ to any open $Y' \subset Y$ with compact
closure has
finite $2$--dimensional Hausdorff measure.
\item $C$ has a positive cohomology assignment.
\end{itemize}
Then the following are true:

\begin{itemize}
\item There is a smooth, complex curve $C^0$ (not necessarily compact) with
a proper,
pseudo-holomorphic map $f\co  C^0 \rightarrow Y$ with $C = f(C^0)$.
\item There is a countable set $\L^0 \subset C^0$ with no accumulation
points such
that $f$ embeds each component $C^0 - \L^0$.
\item Here is an alternate description of the cohomology assignment for
$C$:  Let $\s
\co  D \rightarrow Y$ be an admissible map, and let $\s'$ be any admissible
perturbation of $\s$ which is transverse to $f$ and which is homotopic to
$\s$ via an
admissible homotopy.  Construct the fibered product $T \equiv \{(x,y) \in D
\x C_1:
\s'(x) = \varphi(y)\}$.  This $T$ is a compact, oriented $0$--manifold, so a
finite
set of signed points; and the cohomology assignment gives $\s$ the sum of
the signs
of the points of $T$.
\end{itemize}
\end{proposition}

{\bf Proof of Proposition \ref{prop7.4}}\qua As remarked at the outset,
the proof
of Proposition $6.1$ in \cite{T1} from the revised version in \cite{T6} can be
brought to bear here with negligible modifications.

Lemma \ref{lem7.3} enables Proposition \ref{prop7.4} to be applied to the
set $C$
from Lemma \ref{lem7.1}.  In particular, one can conclude that $C$ is the
image of a
smooth, complex curve $C^0$ via a proper, pseudo-holomorphic
\begin{equation}
\label{eq7.9}
f\co  C^0 \rightarrow X - Z.
\end{equation}
Moreover, $f$ can be taken to be an embedding upon restriction to each
component of
the complement in $C^0$ of a countable set with no accumulation points.  In
particular, it follows that $C$ restricts to any open subset with compact
closure in
$X-Z$ as a pseudo-holomorphic subvariety.

\bs{\bf e) The energy of $C$}

The purpose of this subsection is to state and prove

\begin{lemma}
\label{lem7.5}
The set $C$ from Lemma \ref{lem7.1} is a finite energy, pseudo-holomorphic
subvariety
in the sense of Definition $1.1$.  Furthermore, there are universal constants
$\zeta_1,\zeta_2$ (independent of the metric) such that
\begin{equation}
\label{eq7.10}
\int_C \o \le \zeta_1e_{\o}(s) + \zeta_2\int_X (|R_g| +
|W_g^+|)|\o|\mbox{dvol}_g.
\end{equation}
In this equation, $e_{\o}(s)$ equals the evaluation on the fundamental
class of $X$
of the cup product of $c_1(L)$ with $[\o]$.  Meanwhile, $R_g$ is the scalar
curvature
for the metric $g$, and $W_g^+$ is the metric's self-dual, Weyl curvature.
Also,
$\mbox{dvol}_g$ is the metric's volume form.
\end{lemma}

{\bf Proof of Lemma \ref{lem7.5}}\qua First of all, define an
equivalence relation
on the components of $C^0$ (from (\ref{eq7.9})) by declaring two components
to be
equivalent if their images via $f$ coincide.  The quotient by this equivalence
relation defines another smooth, complex curve, $C_0$, together with a proper,
pseudo-holomorphic map $\varphi\co  C_0 \rightarrow X-Z$ whose image is $C$.
Moreover,
there is a countable set $\L_0 \subset C_0$ which has no accumulation
points and
whose complement is embedded by $\varphi$.

With the preceding understood, it remains only to establish that $C$ has finite
energy.  For this purpose, fix $\d > 0$ and re-introduce the bump function
$\chi_{\d}$.  (Remember that $\chi_{\d}$ vanishes where the distance to $Z$ is
greater than $2\d$, and it equals $1$ where the distance to $Z$ is less
than $\d$.)
Since $\o$ restricts to $C$ as a positive form, it follows that $C$ has
finite energy
if and only if
\begin{equation}
\label{eq7.11}
\lim_{\d \rightarrow 0} \int_C (1-\chi_{\d})\o
\end{equation}
exists, and this limit exists if and only if the set $\left\{ \int_C
(1-\chi_{\d})\o\right\}_{\d > 0}$ is bounded in $[0,\i)$.  Thus, the task
is to find
a $\d$--independent upper bound for this set.

With the preceding understood, remark first that Proposition \ref{prop7.4}
as applied
to $C$ yields
\begin{equation}
\label{eq7.12}
\int_C(1-\chi_{\d})\o \le \int_{C^0} f^*((1-\chi_{\d})\o) = {\cal
F}((1-\chi_{\d})\o).
\end{equation}
Now, given $\e > 0$, the right-hand expression in (\ref{eq7.12}) is no
greater than
\begin{equation}
\label{eq7.13}
{\cal F}_n((1-\chi_{\d})\o) + \e
\end{equation}
when $n$ is sufficiently large.  Moreover, as ${\cal F}_n$ is defined on
all smooth
forms, the first term in (\ref{eq7.13}) is equal to
\begin{equation}
\label{eq7.14}
(4\pi)^{-1}\int_X (1-\chi_{\d})iF_A \wedge \o - (4\pi)^{-1} \int_X (1 -
\chi_{\d})iF_{A^0} \wedge \o.
\end{equation}
with $A = A_n$.

Both terms in (\ref{eq7.14}) are $\d$ dependent and so must be analyzed
further.  In
particular, Lemma \ref{lem7.5} follows with the exhibition of a bound on
these terms
by $\zeta(e_{\o}(s) + \int_X (|R_g| + |W_g^+|)|\o|\mbox{dvol}_g$.

For the right most term, remark that only $F_{A^0}^+$ appears in (\ref{eq7.14})
(since $\o$ is self dual), and a calculation (which is left to the reader)
finds a
universal constant $\zeta$ with the property that
\begin{equation}
\label{eq7.15}
|F_{A^0}^+| \le \zeta(|\nabla\underline{\o}|^2 + |R_g| + |W_g^+|).
\end{equation}
Here, $\underline{\o} \equiv \o/|\o|$.  Equation (\ref{eq7.15}) implies
that the
right most term in (\ref{eq7.14}) is not greater than
\begin{equation}
\label{eq7.16}
\zeta\left(\int_X |\o||\nabla\underline{\o}|^2 +\int_X (|R_g| +
|W_g^+|)|\o|\mbox{dvol}_g\right).
\end{equation}
Here the constant $\zeta$ is also metric independent.  Meanwhile, the first
integral
in (\ref{eq7.16}) is bounded by $\zeta \int_X (|R_g| +
|W_g^+|)|\o|\mbox{dvol}_g$
which can be seen by integrating both sides of (\ref{eq3.7}) and then
integrating by
parts to eliminate the $d^*d|\o|$ term.  This means in particular that the
right most
term in (\ref{eq7.16}) is bounded by a universal multiple of $\int_X (|R_g| +
|W_g^+|)|\o|\mbox{dvol}_g$.

Now consider the left most term in (\ref{eq7.14}).  For this purpose, use
(2.9) to identify the latter with
\begin{equation}
\label{eq7.17}
(8\sqrt{2}\pi)^{-1}r \int_X (1-\chi_{\d})|\o|(2^{-1/2}|\o| - |\a|^2 + |\b|^2).
\end{equation}
One should compare this expression for that implied by (2.9) for
$e_{\o}(s)$,
namely
\begin{equation}
\label{eq7.18}
2^{-1}e_{\o}(s) = (8\sqrt{2}\pi)^{-1}r \int_X |\o|(2^{-1/2}|\o| - |\a|^2 +
|\b|^2).
\end{equation}
In particular, note that (\ref{eq7.18}) implies the identity
\begin{equation}
\label{eq7.19}
(8\sqrt{2}\pi)^{-1}r \int_X |\o||\b|^2 = 4^{-1} e_{\o}(s) -
(16\sqrt{2}\pi)^{-1}r
\int_X |\o|(2^{-1/2}|\o| - |\psi|^2).
\end{equation}
It then follows from the second line of (\ref{eq3.2}) (and the integration of
(\ref{eq3.7})) that
\begin{equation}
\label{eq7.20}
(8\sqrt{2}\pi)^{-1}r \int_X |\o||\b|^2 = \zeta\left(e_{\o}(s) + \int_X (|R_g| +
|W_g^+|)|\o|\mbox{dvol}_g\right).
\end{equation}

This last bound can now be plugged back into bounding (\ref{eq7.17}) since the
expression in (\ref{eq7.17}) is no greater than
\begin{equation}
\label{eq7.21}
(8\sqrt{2}\pi)^{-1}\left(r \int_X \chi_{\d}|\o||(2^{-1/2}|\o| - |\psi|^2)| + 2r
\int_X |\o||\b|^2\right).
\end{equation}
Thus, (\ref{eq7.20}) and the second line of (\ref{eq3.2}) bound the left
most term in
(\ref{eq7.14}) by the required $\zeta\left(e_{\o}(s) + \int_X (|R_g| +
|W_g^+|)|\o|\mbox{dvol}_g\right)$.

\bs{\bf f) Intersections with the linking $2$--spheres}

After Lemma \ref{lem7.5}, all that remains to prove Theorem \ref{th2.2} is to
establish that Lemma \ref{lem7.5}'s pseudo-holomorphic subvariety $C$ has
intersection number equal to $1$ with any $2$--sphere in $X-Z$ which has linking
number $1$ with $Z$.

To prove this last assertion, consider that each $\a_n$ defines a section
of the
bundle $E$ whose square is given in (\ref{eq6.1}) as $L|_{X-Z}K$.  Now, $L$ is
trivial over a linking $2$--sphere for $Z$ as $L$ is a bundle over the whole
of $X$.
This means that $E$ restricts to a linking $2$--sphere of $Z$ as the square
root of
the restriction of $K$.  Moreover, according to Lemma \ref{lem2.1}, the
restriction of
$K$ to a linking $2$--sphere has degree $2$, so $E$ restricts to such a
$2$--sphere
with degree $1$.  This means that the current ${\cal F}_n$ in (\ref{eq7.1})
evaluates
as $1$ on any closed form which represents the Thom class of a fixed, linking
$2$--sphere.  In particular, the same must be true for a limit current
${\cal F}$, and
the nature of the convergence in Proposition \ref{prop7.4} implies that $C$ has
intersection number $1$ with such a $2$--sphere.

\rk{Achnowledgement}This research was supported in part by the National
Science Foundation.

\end{document}